\documentclass[12pt,a4paper]{article}
\usepackage{amssymb,amsmath,bm}
\usepackage{fullpage}
\usepackage{mathtools}
\usepackage{booktabs}

\usepackage{algorithm}
\usepackage{algpseudocode}
\usepackage{natbib}

\usepackage{graphicx}
\usepackage{tikz}
\usepackage{pgfplots}

\newcommand{\R}{\ensuremath{\mathbb{R}}}

\usepackage{hyperref}
\hypersetup{
    colorlinks=true,
    citecolor=black,
    urlcolor=blue,
    linkcolor=black,
    pdfborder={0 0 0},
}

\usepackage{environ,tabularx}
\makeatletter
\newcommand{\problemtitle}[1]{\gdef\@problemtitle{#1}}
\newcommand{\probleminput}[1]{\gdef\@probleminput{#1}}
\newcommand{\problemquestion}[1]{\gdef\@problemquestion{#1}}
\NewEnviron{problem}{
	\problemtitle{}\probleminput{}\problemquestion{}
	\BODY
	\noindent
	\par\addvspace{.5\baselineskip}
	\noindent
	\begin{tabularx}{0.8\textwidth}{@{\hspace{\parindent}} l X c}
		\multicolumn{2}{@{\hspace{\parindent}}l}{\sc \@problemtitle} \\
		\textbf{Input:} & \@probleminput \\
		\textbf{Question:} & \@problemquestion
	\end{tabularx}
	\par\addvspace{.5\baselineskip}
}
\makeatother

\begin{document}
\title{Interior point methods can exploit structure of convex piecewise linear functions with application in radiation therapy}
\author{Bram L. Gorissen}
\date{}
\maketitle
\abstract{Auxiliary variables are often used to model a convex piecewise linear function in the framework of linear optimization. This work shows that such variables yield a block diagonal plus low rank structure in the reduced KKT system of the dual problem. We show how the structure can be detected efficiently, and derive the linear algebra formulas for an interior point method which exploit such structure. The structure is detected in 36\% of the cases in Netlib. Numerical results on the inverse planning problem in radiation therapy show an order of magnitude speed-up compared to the state-of-the-art interior point solver CPLEX, and considerable improvements in dose distribution compared to current algorithms.}

\begin{tikzpicture}[remember picture,overlay]
\node[anchor=south,yshift=10pt] at (current page.south) {\fbox{\parbox{\dimexpr\textwidth-\fboxsep-\fboxrule\relax}{\footnotesize This is an author-created, un-copyedited version of an article accepted for publication in SIAM Journal on Computing  \href{https://doi.org/10.1137/21M1402364}{DOI:10.1137/21M1402364}. The appendix and the Python code that accompany the manuscript on arXiv are not part of the publication.}}};
\end{tikzpicture}
\section{Introduction}

Linear optimization is a versatile modeling framework with widespread applications. Its success is at least partially explained by the availability of well performing and stable algorithms, such as the simplex method and interior point methods (IPMs), both of which have seen tremendous improvements over the past decades \citep{bixby2012brief}.

Formulating a linear optimization problem often follows certain patterns. A set of standard reformulations can be used to linearize expressions, such as convex piecewise linear constraints, which adds auxiliary variables. Due to such reformulations, linear optimization problems often have a specific structure that is currently ignored. Since models are generated with a modeling environment such as AMPL \citep{fourer1997ampl} or by hand, constraints have a natural order, which allows the structure to be detected in a greedy fashion.

In each iteration of an IPM that solves a linear optimization problem with coefficient matrix $\bm{A}$, a linear system has to be solved of which the sparsity structure is the same as that of $\bm{A}\bm{A}^T$. The solution time of an IPM is therefore closely related to the structure of $\bm{A}\bm{A}^T$. The structure we shall exploit is block diagonal plus low rank, except for the first block row and column:
\begin{align}
\bm{A}\bm{A}^T = \begin{pmatrix}
(\bm{AA}^T)_{11}   & (\bm{AA}^T)_{12}    & \hdots & (\bm{AA}^T)_{1K} \\
(\bm{AA}^T)_{12}^T & \bm{D}_2 + \bm{R}_2 &        & \bm{O} \\
\vdots             &                     & \ddots & \\
(\bm{AA}^T)_{1K}^T & \bm{O}              &        & \bm{D}_K + \bm{R}_K
\end{pmatrix}, \label{coeff:structure}
\end{align}
where $\bm{D}_k$ are diagonal and $\bm{R}_k$ are low rank matrices for $k=2,\ldots,K$. This structure has not been recognized before to the best of our knowledge, although a pure block diagonal plus low rank system has appeared in specific problems, e.g., multiclass support vector machines \citep[\S 1.3.3]{andersen2011interior} and large-scale power systems \citep{zhang2017robust}. This structure is a generalization of the master problem in the decomposition by \citet{gondzio1997using}. Similarities between our approach and other IPMs are described in Section \ref{sec:comparison}.

We perform an efficient row reduction of $\bm{A}\bm{A}^T$ to obtain a smaller system where the auxiliary variables are eliminated. Let $m_k$ be the size of the $k$-th diagonal block matrix and let $m = \sum_{k=1}^K m_k$. Without exploiting structure, solving a linear system with coefficient matrix \eqref{coeff:structure} takes $\mathcal{O}(m^3)$ steps. We consider two extremes. Suppose $m_2=m_3=\ldots=m_K=1$, then nothing can be gained. If on the other hand $K=2$, $m_1=1$, $m_2=m-1$ and the rank of $\bm{R}_2$ is $p_2$, the system can be solved in $\mathcal{O}(p_2^3 m)$ steps by reducing the coefficient matrix to the scalar $(\bm{AA}^T)_{11} + (\bm{AA}^T)_{12} (\bm{D}_2^2 + \bm{R}_2)^{-1}(\bm{AA}^T)_{12}^T$. So, if $p_2$ is small, the solution method scales linearly in the number of auxiliary variables.

An important application of this method is radiation therapy, which is a commonly applied treatment for cancer. A challenge in day-to-day clinical care is to get ionizing radiation into the tumor while minimizing the exposure of surrounding organs. Due to the large number of degrees of freedom, optimization models have been used for decades to design a treatment plan based on dosimetric endpoints. Each relevant organ is discretized into cubes (voxels), indexed by $i$. The dose in a voxel is linear in the fluence vector $\bm{x} \in \mathbb{R}^n_+$ for most treatment modalities \citep{bortfeld1990methods,lessard2001inverse,lomax1999intensity}:
\begin{align}\label{eq:dose}
d_i = \sum_j D_{ij} x_j.
\end{align} 
The dose kernel matrix $\bm{D}$ is a parameter that models the relationship between the dose in voxel $i$ and the fluence from position $j$. Dose is measured in Gray (Gy). Dosimetric endpoints are formulated in the space of $d_i$. Consider the following objective functions \citep{romeijn2006new}:
\begin{enumerate}
\item Mean underdose: $\sum_i \max\{0, p - d_i\} / n$, where $i$ sums over all $n$ voxels in the tumor and $p$ denotes the prescribed tumor dose.
\item Conditional value at risk (CVaR) at level $\alpha$: $\min_y \{ y + \sum_i \max\{0, y - d_i \} / ((1-\alpha)n) \}$, which is the mean dose in the $\alpha\%$ region of the tumor that receives the least dose.
\item Dose-volume histogram (DVH) statistic at level $\alpha$: this is the maximum dose in the $\alpha\%$ region of the tumor that receives the least dose. This can be modeled with binary variables or solved heuristically \citep{mukherjee2019integrating}, and both methods give rise to the sum of $\max\{\cdot\}$ expressions.
\end{enumerate}
While these objectives are dose promoting, similar objectives with $\max\{\cdot\}$ expressions can be defined for the organs at risk (OARs), where the goal is to minimize dose. In summary, whether these three dose based functions are used for the target or for OARs, they use auxiliary variables to reformulate a $\max\{\cdot \}$ expression, and give rise to exploitable structure when they are used in an objective or constraint.

Commercially available treatment planning systems (TPSs) often use L-BFGS, the conjugate gradient method or simulated annealing, and rely on techniques such as squaring variables to ensure nonnegativity or using log-sum-exp to combine many linear constraints into one convex constraint \citep{fredriksson2014critical,karabis2009optimization}. Although these algorithms are fast, they are known to produce suboptimal results. It is likely that the potential of radiation therapy is not fully utilized at least for some patients due to limitations of optimization algorithms.

There are three approaches to address the multiobjective nature of treatment planning: 1.~construct a Pareto frontier and allow the user to navigate over this frontier, 2.~use lexicographic optimization, or 3. let the user manually tune weights for different objectives \cite{breedveld2019multi}. These approaches are most effective when there is a limited number of meaningful objectives. The CVaR objective was proposed for radiation therapy already in 2003 \citep{romeijn2003novel}, but is not supported by any TPS. Recently, CVaR was shown to make the trade-off between objectives easier and to yield better treatment plans than other objectives \citep{engberg2017explicit}.

An implementation of the interior point method proposed in this paper named Nymph was incorporated in the TPS Astroid (.decimal, inc., Sanford, FL, USA) and is used at the Francis H.~Burr Proton Therapy Center (Massachusetts General Hospital, Boston, MA, USA), where it has fully replaced the previous algorithm ART3+O \citep{chen2010fast}, which is described in Section \ref{sec:radiation}. Due to support for CVaR, for mean underdose constraints, and an optimality guarantee, creating a treatment plan takes less effort than before. To the best of our knowledge, it is the first clinical application of CVaR. A fully functional and unrestricted standalone version of Nymph is freely available for research purposes from the author's website \citep{Gorissen2020}.

We use the following notation. Matrices are shown in bold. The elements of the matrix $\bm{O}$ and the vector $\bm{0}$ are all zero. The matrix $\bm{I}_k$ is an identity matrix of dimension $k$. The vector $\bm{1}_k$ denotes a $k$-dimensional vector whose elements are all $1$. The operator $\otimes$ denotes the Kronecker product.

The remainder of the paper is organized as follows. Section \ref{sec:methods} contains the main results. Section \ref{sec:structure} shows that the structure in \eqref{coeff:structure} appears from convex piecewise linear constraints, and provides examples of how such constraints arise. It is argued that the computationally most appealing cases are constraints on the sum of absolute values (the $\ell_1$ norm) or the sum of terms that each is the maximum of 0 and a linear function of the decision variables. Section \ref{sec:detection} proves that optimally detecting the structure in a given linear optimization problem is NP-hard, but provides a procedure for detecting the structure in a greedy way. Section \ref{sec:fullIPM} gives an explicit derivation of the linear system that determines the search direction for an interior point method. Section \ref{sec:comparison} provides a comparison with other IPMs that exploit structure. Section \ref{sec:examples} shows an evaluation of the greedy detection algorithm on the Netlib/lp test set. It also demonstrates the efficacy of the proposed method for radiation therapy optimization based on a quantitative comparison with CPLEX.
\section{Methods}\label{sec:methods}
\subsection{Problem structure}\label{sec:structure}
Consider a linear optimization problem in standard form and its dual:
\begin{align*}
\min_{\bm{x} \in \R^n} \quad & \bm{c}^T \bm{x} & \max_{\bm{y} \in \R^m, \bm{s} \in \R^n} \quad & \bm{b}^T \bm{y} \\
\mathrm{s.t. } \quad & \bm{Ax} = \bm{b} & \mathrm{s.t.} \quad & \bm{A}^T\bm{y}+\bm{s}=\bm{c} \\
&\bm{x} \geq \bm{0} & & \bm{s}\geq\bm{0},
\end{align*}
where $\bm{c} \in \mathbb{R}^{n}$ defines the objective function, and $\bm{A} \in \mathbb{R}^{m \times n}$ and $\bm{b} \in \mathbb{R}^{m}$ model the constraints in the primal. The structure in \eqref{coeff:structure} arises if $\bm{A}$ is of the form:
\begin{align}
\bm{A}=\begin{pmatrix}
\bm{A}_{11} & \bm{A}_{12} & \hdots & \bm{A}_{1K} \\
\bm{A}_{21} & \bm{E}_2    &        & \bm{O} \\
\vdots      &             & \ddots & \\
\bm{A}_{K1} & \bm{O}      &        & \bm{E}_K
\end{pmatrix}, \label{coeff:structureA}
\end{align}
where the matrices $\bm{A}_{k1} \in \mathbb{R}^{m_k \times n_1}$ ($k=2,\ldots,K$) are structurally mutually orthogonal ($\bm{A}_{p1} \bm{W} \bm{A}_{q1}^T = \bm{O}$ for all $p,q \geq 2$ with $p\neq q$, and all diagonal matrices $\bm{W}$) and have a limited number of nonzero columns, and $\bm{E}_k \in \mathbb{R}^{m_k \times n_k}$ is such that $\bm{E}_k\bm{E}_k^T$ is diagonal. There is no restriction on the matrices $\bm{A}_{1k} \in \mathbb{R}^{m_1 \times n_k}$ ($k=2,\ldots,K$).

This structure of $\bm{A}$ occurs naturally in many linear optimization problems. We conjecture that the most common occurrence is due to the linear reformulation of convex piecewise linear (CPL) constraints with coefficient matrices $\bm{F}^j$:
\begin{align}
\sum_{i=1}^P \max_{1\leq l \leq L} \left\{\left(\bm{F}^j_{il}\right)^T \bm{y}\right\}	 \leq c_j \; \forall j \quad \Leftrightarrow \quad \sum_{i=1}^P z_{ij} \leq c_j, \; \left(\bm{F}^j_{il}\right)^T\bm{y} \leq z_{ij} \; \forall i,j,l. \label{eq:piecewise}
\end{align}
The number of blocks $K$ in \eqref{coeff:structureA} is equal to the number of reformulated CPL constraints. For example, if an optimization problem has one CPL constraint, we can drop the index $j$ and constraint \eqref{eq:piecewise} can be expressed as a set of linear constraints by adding slack variables $\bm{s}$:
\begin{align}
\setcounter{MaxMatrixCols}{20}
\begin{pmatrix}
\bm{0}  & \bm{F}_{11} & \cdots & \bm{F}_{1L} & \bm{F}_{21} & \cdots & \bm{F}_{2L} & & \bm{F}_{P1} & \cdots & \bm{F}_{PL} \\
1 & -1 & \cdots & -1 \\
1 &   &        &   & -1 & \cdots & -1 \\
1 &   &        &   &   &        &   & \ddots \\
1 &   &        &   &   &        &   &        & -1 & \cdots & -1
\end{pmatrix}^T
\begin{pmatrix} \bm{y} \\ \bm{z} \end{pmatrix} + \bm{s} = \begin{pmatrix}\bm{c} \\ \bm{0}\end{pmatrix}, \label{eq:exampleA}
\end{align}
which is exactly the structure of the dual problem where $\bm{y}$ is substituted by $(\bm{y},\bm{z})$ and $\bm{c}$ is substituted by $(\bm{c},\bm{0})$. Indeed $\bm{A}_{21}$ has only one nonzero column, and $\bm{E}_2 = -\bm{I}_P \otimes \bm{1}_L^T$ is such that $\bm{E}_2\bm{E}_2^T$ is diagonal.

CPL constraints arise naturally in many optimization problems:
\begin{enumerate} \setlength\itemsep{0em}
	\item to model cost functions, such as for holding and backlogging costs in inventory problems:
	$$\sum_{t=0}^T \max\{(x_0 + \sum_{i=0}^t x_i - d_i)c_h, -(x_0 + \sum_{i=0}^t x_i - d_i)c_b \},$$
	\item to model absolute values, such as for $\ell_1$ regression \citep{charnes1955optimal}:
	$$||\bm{Ax}-\bm{b}||_1,$$
	\item to model conditional value at risk \citep{rockafellar2000optimization}:
	$$\alpha + \frac{1}{1-\beta} \sum_i \max\{0, \bm{f}_i^T\bm{x}-\alpha\}, $$
	\item to model soft constraints, such as in radiation therapy \citep{shepard1999optimizing}:
	$$\sum_i \max\left\{ 0, \alpha_i \left(\sum_j d_{ij}x_j - L_i\right), \beta_i \left(U_i - \sum_j d_{ij}x_j\right) \right\}\text{, and} $$
	\item to model the safety factor in Robust Optimization under box uncertainty or budget uncertainty \citep{bertsimas2004price}:
	$$(\bm{a}+\bm{A}\bm{\zeta})^T\bm{x} \leq b \; \forall \bm{\zeta} : ||\bm{\zeta}||_1 \leq \Gamma, ||\bm{\zeta}||_\infty \leq \rho \Leftrightarrow \bm{a}^T\bm{x} +\Gamma || \bm{A}^T\bm{x} - \bm{z}||_\infty + \rho ||\bm{z}||_1 \leq b.$$
\end{enumerate}

\subsection{Structure detection}\label{sec:detection}
The structure can be detected at three levels. We discuss these levels top-down, starting with the user level. Since the desired structure \eqref{coeff:structureA} is often the result of reformulating a CPL constraint \eqref{eq:piecewise}, it is reasonable to request the modeler to pass that information to the solver. The solver then knows explicitly which blocks can be eliminated. However, it would be wasteful to solely rely on the modeler to provide structural information. They may not be aware that such information can be useful, and existing models and software may not receive performance updates.

At the intermediate level, the structure can be detected in an algebraic modeling environment. Many optimization problems are entered into a computer in an algebraic way with software such as AMPL and AIMMS \citep{fourer1997ampl}. This software can access constraints in their algebraic form, and  can therefore detect the structure in \eqref{eq:piecewise} by looking for constraints that have: (a) an indexed (auxiliary) variable on one side of an inequality, (b) a summation over one or more indices on the other side of the inequality, and (c) a separate constraint where the auxiliary variables are summed. This structure can then be passed to the solver via a special interface. This approach was proposed by \cite{fragniere2000structure} for a block-angular structure, but unfortunately never gained traction.

At the bottom level it is therefore useful to have a procedure that can efficiently detect exploitable structure directly from the coefficent matrix. Detecting the structure in an optimal way should take into account that the rows and columns of $\bm{A}$ may be permuted, which makes this approach infeasible in practice. Even detecting the largest single block when $K=2$ and $\bm{A}_{21} \in \mathbb{R}^{m_1 \times 1}$ is NP-hard. Consider the following two decision problems:
\begin{problem}
	\problemtitle{IndependentSet}
	\probleminput{An undirected graph $G = (V,E)$ and an integer $k\in\mathbb{N}$.}
	\problemquestion{Is there a subset $S \subseteq V$ such that $|S| \geq k$ and $\{v_1,v_2\} \not \in E$ for all $v_1,v_2 \in S$?}
\end{problem}
\begin{problem}
	\problemtitle{BlockStructure}
	\probleminput{A matrix $\bm{B}\in\mathbb{R}^{m\times n}$ and an integer $k\in\mathbb{N}$.}
	\problemquestion{Is there a subset of the rows $I \subseteq \{1,2,\ldots,m\}$ and a column $q \in \{1,2,\ldots,n\}$ such that $|I| \geq k$ and
	
	$|\{ i \in I : \bm{B}_{ij} \neq 0\}| \leq 1$ for all $j \neq q$?}
\end{problem}
The latter problem can identify the structure of \eqref{coeff:structureA} for $K=2$: $I$ defines the rows of the second block row (which may require permuting rows such that the rows in $I$ are at the bottom), and $q$ defines the first block column (which may require permuting columns). Indeed the structure will ensure that $\bm{A}_{21}$ is rank 1 and that $\bm{E}_k\bm{E}_k^T$ is diagonal.

Note that {\sc BlockStructure} is in NP. It is well known that {\sc IndependentSet} is NP-complete, since it is equivalent to {\sc Clique} on the complementary graph and {\sc Clique} is NP-complete \citep{karp1972reducibility}. We provide a polynomial-time reduction from {\sc IndependentSet} to {\sc BlockStructure}, ignoring the trivial case $k\leq 1$: take $m=|V|$, $n=|E|+1$, $\bm{A}_{i1}=1$ for all $i$, and $\bm{A}_{ij}=1$ $(j\geq 2$) if vertex $i$ is in edge $j-1$. The correctness of the reduction follows trivially from the fact that the only $q$ that can be chosen is $q=1$, and that $\{v_1,v_2\} \in E$ iff $\bm{A}_{v_1,j} = \bm{A}_{v_2,j} = 1$ for some $j\in\{2,\ldots,|E|+1\}$. {\sc BlockStructure} is therefore NP-complete.

Although it may be impossible to identify an optimal structure, we rarely face random optimization problems. Instances are often derived from a structural model, either on paper or in a modeling environment such as AMPL. Using the structural model to create the constraint matrix yields a natural ordering of the rows. In particular, it is likely that the linear constraints that model a CPL constraint appear consecutively. We therefore propose to go over the rows of $\bm{A}$ in a single pass and test whether consecutive constraints belong to the same block matrix. The algorithm does not assume a particular ordering of the constraints within each block, and the algorithm is invariant under permutations of the columns of $\bm{A}$. The steps of our method are shown as Alg.~\ref{alg:structure}. The remainder of this section explains this algorithm.

The algorithm has two parameters. The parameter $m^{\min} \geq 2$ is a filter on the minimum number of rows of a block. Blocks with fewer rows are considered part of $\bm{A}_{11}$. Small blocks offer little computational benefit, and explicitly not detecting them makes it more likely to find other blocks. The parameter $J^{\max}$ is the maximum number of nonzero elements in a row of $\bm{A}$ to be considered part of a block. A small number of nonzero elements is an indication that $\bm{A}_{k1}$ has low rank and that $\bm{E}_k^T \bm{E}_k$ has a limited number of off-diagonal entries. The algorithm returns the row indices of the blocks ($\mathcal{I}_k$), the column indices of the blocks on the diagonal ($\mathcal{E}_k$) and the column indices of the structurally mutually orthogonal matrices ($\mathcal{M}_k$). So, the rows and columns of $\bm{A}_{21}$ are $\mathcal{I}_2$ and $\mathcal{M}_2$ whereas the columns of $\bm{A}_{12}$ and $\bm{E}_2$ are given by $\mathcal{E}_2$. The row and column indices of $\bm{A}_{11}$ are the complements of $\cup_{k=2}^K \mathcal{I}_k$ and $\cup_{k=2}^K \mathcal{E}_k \cup \mathcal{M}_k$, respectively. For example, for the coefficient matrix in \eqref{eq:exampleA}, the algorithm returns $K=2$, $\mathcal{I}_2=\{m+1,m+2,\ldots,m+P\}$, $\mathcal{E}_2=\{2,3,\ldots,PL\}$, and $\mathcal{M}_2=\{1\}$, where $m$ denotes the dimension of $\bm{y}$.

The algorithm starts with an attempt to identify the first block ($k=2$) and initializes the index sets as empty sets. The main loop goes over the rows, starting with the second row. If the number of nonzero elements in the current and previous row are both less than $J^{\max}$ and if the previous row is not part of an existing block, those two rows form the start of block $k$ where the columns for $\bm{A}_{k1}$ are those for which there is a nonzero entry in both rows. For the next row of $\bm{A}$ we check if it fits the pattern, i.e., if it also has a nonzero. If we find a row that does not fit in the pattern or has too many nonzero entries, we either create a new block if the current block has at least $m^{\min}$ rows and the columns are not part of a previously detected block, or erase the current attempt to create a block and start over.

The algorithm is denoted in set notation. One can make further assumptions about the order of the auxiliary variables that model a CPL constraint \eqref{eq:piecewise}. This allows for a more efficient implementation than with basic set operations. The most significant improvement is storing only the smallest and largest element of the set $\mathcal{M}_k$, and assuming all elements in between are in $\mathcal{M}_k$. Determining if a set is a subset of another set (within line 13 of the algorithm) then only requires two comparisons between numbers.

\begin{algorithm}
\caption{Detecting structure}\label{alg:structure}
\begin{algorithmic}[1]
\State Input: $\bm{A}\in\mathbb{R}^{m \times n}$
\State Output: $K \in \mathbb{N}$ and $\{ (\mathcal{I}_k,\mathcal{E}_k,\mathcal{M}_k) : k \in \{2,\ldots,K\} \}$
\State $\{j : \bm{A}_{ij} \neq 0 \} \to \mathcal{J}_i$ for $i=1,\ldots,m$ (nonzero columns in row $i$)
\State $2 \to k$ (index of current block)
\State $\emptyset \to \mathcal{D}$ (set of columns that are part of a detected block)
\State $\emptyset \to \mathcal{I}_k$, $\emptyset \to \mathcal{E}_k$, and $\emptyset \to \mathcal{M}_k$ (row and column indexes of first block)
\State $\emptyset \to \mathcal{I}_1$ (inititalize $\mathcal{I}_{k-1}$)
\For{$i=2$ to $m$}
\If{$|\mathcal{J}_i| \leq J^{\max}$ \text{\bf and } $|\mathcal{J}_{i-1}| \leq J^{\max}$ \text{\bf and } $i-1\not\in \mathcal{I}_{k-1}$}
\If{$\mathcal{I}_{k} = \emptyset$}
\State $\{i-1,i \} \to \mathcal{I}_k, \mathcal{J}_{i-1} \cap \mathcal{J}_i \to \mathcal{M}_k, (\mathcal{J}_i \cup \mathcal{J}_{i-1}) \backslash \mathcal{M}_k \to \mathcal{E}_k$
\Else
\If{$\mathcal{J}_{i-1} \cap \mathcal{J}_i \subseteq \mathcal{M}_k$ \text{\bf and } $\mathcal{E}_k \cap (\mathcal{J}_i \backslash \mathcal{M}_k) = \emptyset$}
\State $\mathcal{I}_k \cup \{i\} \to \mathcal{I}_k, \mathcal{M}_k \cup (\mathcal{J}_i \backslash \mathcal{M}_k) \to \mathcal{E}_k$
\Else
\State {\bf goto} line 20
\EndIf
\EndIf
\Else
\If{$|\mathcal{I}_k| \geq m^{\min}$ {\bf and } $(\mathcal{E}_k \cup \mathcal{M}_k)\cap \mathcal{D}=\emptyset$}
\State $\mathcal{D} \cup \mathcal{E}_k\cup \mathcal{M}_k\to \mathcal{D}$
\State $k+1 \to k$
\EndIf
\State $\emptyset \to \mathcal{I}_k$, $\emptyset \to \mathcal{E}_k$, and $\emptyset \to \mathcal{M}_k$
\EndIf
\EndFor
\State {\bf if} $|\mathcal{I}_k| \geq m^{\min}$ {\bf and } $(\mathcal{E}_k \cup \mathcal{M}_k)\cap \mathcal{D}=\emptyset$ {\bf then} $k \to K$ {\bf else}  $k-1 \to K$
\end{algorithmic}
\end{algorithm}

\subsection{Interior point method}\label{sec:fullIPM}
We assume that the coefficient matrix $\bm{A}$ has the structure \eqref{coeff:structureA}. It may be necessary to dualize the problem to get it in the required format. The logarithmic barrier formulation with parameter $\mu$, justifying the name {\it interior} point method, is given by:
\begin{align*}
\min_{\bm{x} \in \R^n} \left\{ \bm{c}^T \bm{x} - \mu \sum_{i=1}^n \log(x_i) \mid
\bm{Ax} = \bm{b} \right\}.
\end{align*}
The KKT conditions, necessary and sufficient for optimality, are:
\begin{align*}
\bm{Ax} = \bm{b}, \quad
\bm{c} - \bm{A}^T \bm{y} - \bm{s} = \bm{0}, \quad
\bm{X}\bm{S}\bm{e} - \mu \bm{e} = \bm{0} \text{, and }
\bm{x},\bm{s} \geq \bm{0},
\end{align*}
where we adopt the notation that $\bm{X}$ and $\bm{S}$ are square diagonal matrices with $\bm{x}$ and $\bm{s}$ on the diagonal, respectively. The essence of a primal-dual IPM is remarkably simple and is displayed as Alg.~\ref{basic_ipm}. The first step is generating a starting point that is sufficiently far from the boundary with, e.g., the heuristic in \citep{mehrotra1992implementation}. The remainder of the algorithm is a loop that continues until the KKT conditions are satisfied with sufficient numerical precision. In each iteration, the weight of the logarithmic barrier $\mu$ is reduced, after which a search direction is determined and a step of length $\alpha$ is taken such that $\bm{x}$ and $\bm{s}$ remain strictly positive. Despite the simple basics, a good implementation is more complex due to dealing with free variables \citep{vanderbei1999loqo}, detecting infeasibility or unboundedness \citep{andersen2000mosek}, factorizing the coefficient matrix \citep{andersen2000mosek,vanderbei1999loqo}, exploiting sparsity \citep{andersen2000mosek,vanderbei1999loqo}, improving the search direction with predictor-corrector methods \citep{colombo2008further,gondzio1996multiple,mehrotra1992implementation}, etc. The basic IPM presented here is sufficient to demonstrate our proposed method, although our method fits seamlessly in more sophisticated IPMs.

\begin{algorithm}
\caption{Basic outline of an IPM}\label{basic_ipm}
\begin{algorithmic}
\State Input: $\bm{A} \in \mathbb{R}^{m \times n}$, $\bm{b} \in \mathbb{R}^{m}$, $\bm{c} \in \mathbb{R}^{n}$
\State Output: $\bm{x}$ that solves $\min_{\bm{x} \in \R^n} \{ \bm{c}^T \bm{x} : \bm{Ax} = \bm{b}, \bm{x} \geq \bm{0} \}$ and certificate $(\bm{y},\bm{s})$
\State Generate starting point $(\bm{x},\bm{y},\bm{s})$
\Repeat
\State Determine barrier parameter $\mu$ based on the centrality of $\bm{x}$ and $\bm{s}$
\State Determine search direction $(\Delta \bm{x},\Delta \bm{y},\Delta\bm{s})$
\State Determine step length $\alpha$ such that $\bm{x}+\alpha\Delta \bm{x} > 0$ and $\bm{s}+\alpha\Delta \bm{s} > 0$
\State $(\bm{x},\bm{y},\bm{s}) \gets (\bm{x},\bm{y},\bm{s}) + \alpha (\Delta \bm{x},\Delta \bm{y},\Delta\bm{s})$
\State $\bm{r}^{primal} \gets \bm{b} - \bm{Ax}$, $\bm{r}^{dual} \gets 
\bm{c} - \bm{A}^T \bm{y} - \bm{s}$, $\bm{r}^{comp} \gets \mu \bm{e} - \bm{X}\bm{S}\bm{e}$
\Until{$||\bm{r}^{primal}||<\varepsilon_p (1+||\bm{b}||) \wedge ||\bm{r}^{dual}||<\varepsilon_d (1+||\bm{c}||) \wedge ||\bm{r}^{comp}||<\varepsilon_c (1+||\bm{c}^T\bm{x}||)$}
\end{algorithmic}
\end{algorithm}

\noindent The search direction is found by linearizing the KKT conditions at the current iterate:
\begin{align*}
\begin{pmatrix}
\bm{A} & \bm{O} & \bm{O} \\
\bm{O} & \bm{A}^T & \bm{I} \\
\bm{S} & \bm{O} & \bm{X}
\end{pmatrix}
\begin{pmatrix}
\Delta \bm{x} \\
\Delta \bm{y} \\
\Delta \bm{s}
\end{pmatrix}
=
\begin{pmatrix}
\bm{b} -\bm{Ax} \\
\bm{c} - \bm{A}^T \bm{y} - \bm{s} \\
\mu \bm{e} -\bm{X}\bm{S}\bm{e} \\
\end{pmatrix}
=
\begin{pmatrix}
\bm{r}^{primal} \\
\bm{r}^{dual} \\
\bm{r}^{comp} \\
\end{pmatrix}.
\end{align*}
This system can be reduced by eliminating $\Delta \bm{s}$ and $\Delta \bm{x}$:
\begin{align*}
\Delta \bm{s} = \bm{r}^{dual} - \bm{A}^T \Delta \bm{y}, \text{ and }
\Delta \bm{x} = \bm{S}^{-1} \bm{r}^{comp} - \bm{X}\bm{S}^{-1} \Delta \bm{s}.
\end{align*}
These reduction steps are computationally cheap to carry out since the inverted matrices are diagonal. The linear system then becomes the normal equation
\begin{align}
&\bm{A} \bm{X}\bm{S}^{-1} \bm{A}^T \Delta \bm{y} = \bm{r}^{primal} + \bm{A} \left( \bm{X}\bm{S}^{-1} \bm{r}^{dual} - \bm{S}^{-1} \bm{r}^{comp} \right). \label{eq:norm}
\end{align}
This system, denoted here as $\bm{N} \Delta \bm{y} = \bm{r}$ for simplicity, can be solved by first taking the Cholesky decomposition $\bm{N} = \bm{L}\bm{L}^T$ where $\bm{L}$ is a lower triangular matrix, then solving $\bm{L} \bm{z} = \bm{r}$, and finally solving $\bm{L}^T \Delta \bm{y} = \bm{z}$. However, due to the assumed structure on $\bm{A}$, we shall first reduce the system before doing a Cholesky factorization. After substituting \eqref{coeff:structureA} into \eqref{eq:norm}, the coefficient matrix becomes:
\begin{align}
&\begin{pmatrix}
\sum_{k=1}^K \bm{A}_{1k} \bm{W}_k \bm{A}_{1k}^T &
\bm{A}_{11}\bm{W}_1\bm{A}_{21}^T+\bm{A}_{12}\bm{W}_2\bm{E}_2^T &
\hdots &
\bm{A}_{11}\bm{W}_1\bm{A}_{K1}^T+\bm{A}_{1K}\bm{W}_K\bm{E}_K^T \\
\bm{A}_{21}\bm{W}_1\bm{A}_{11}^T + \bm{E}_2 \bm{W}_2 \bm{A}_{12}^T &
\bm{A}_{21}\bm{W}_1\bm{A}_{21}^T + \bm{E}_2 \bm{W}_2 \bm{E}_2^T &
&
\bm{O} \\
\vdots & & \ddots \\
\bm{A}_{K1}\bm{W}_1\bm{A}_{11}^T + \bm{E}_K \bm{W}_K \bm{A}_{1K}^T &
\bm{O} &
&
\bm{A}_{K1}\bm{W}_1\bm{A}_{K1}^T + \bm{E}_K \bm{W}_K \bm{E}_K^T &
\end{pmatrix}, \label{eq:norm:coeff}
\end{align}
where $\bm{W} = \bm{X}\bm{S}^{-1}$ for brevity, with block matrices $\bm{W}_1, \ldots, \bm{W}_K$ on the diagonal. The variables $\Delta\bm{y}_k$ ($k=2,\ldots,K$) can be eliminated via:
\begin{align}
\Delta\bm{y}_k = \left( \bm{A}_{k1}\bm{W}_1\bm{A}_{k1}^T + \bm{E}_k \bm{W}_k \bm{E}_k^T \right)^{-1} \left( \bm{r}_k - \left(\bm{A}_{k1}\bm{W}_1\bm{A}_{11}^T + \bm{E}_k \bm{W}_k \bm{A}_{1k}^T\right) \Delta\bm{y}_1 \right). \label{eq:deltay}
\end{align}
Defining $\bm{B}_k = \bm{A}_{11}\bm{W}_1\bm{A}_{k1}^T+\bm{A}_{1k}\bm{W}_k\bm{E}_k^T$, the reduced coefficient matrix is:
\begin{align}
\sum_{k=1}^K \bm{A}_{1k} \bm{W}_k \bm{A}_{1k}^T - 
\sum_{k=2}^K \bm{B}_k \left( \bm{A}_{k1}\bm{W}_1\bm{A}_{k1}^T + \bm{E}_k \bm{W}_k \bm{E}_k^T \right)^{-1} \bm{B}_k^T. \label{coeff:reduced:beforesmw}
\end{align}
Let $p_k$ denote the number of nonzero columns of $\bm{A}_{k1}\bm{W}_1\bm{A}_{k1}^T$ which was assumed to be low. There is a $\bm{V}_k \in \mathbb{R}^{m_k \times p_k}$ such that $\bm{A}_{k1}\bm{W}_1\bm{A}_{k1}^T = \bm{V}_k\bm{V}_k^T$. The Sherman-Morrison-Woodbury formula \citep{guttman1946enlargement} allows us to express the reduced coefficient matrix \eqref{coeff:reduced:beforesmw} as:
\begin{align}
& \sum_{k=1}^K \bm{A}_{1k} \bm{W}_k \bm{A}_{1k}^T - \sum_{k=2}^K \bm{B}_k \left(  \bm{E}_k \bm{W}_k \bm{E}_k^T \right)^{-1} \bm{B}_k^T \notag \\
& \quad + \sum_{k=2}^K \bm{B}_k \left(  \bm{E}_k \bm{W}_k \bm{E}_k^T \right)^{-1}
\bm{V}_k
\left(\bm{I} + \bm{V}_k^T \left(  \bm{E}_k \bm{W}_k \bm{E}_k^T \right)^{-1} \bm{V}_k\right)^{-1}
\bm{V}_k^T
\left(  \bm{E}_k \bm{W}_k \bm{E}_k^T \right)^{-1} \bm{B}_k^T. \label{coeff:reduced}
\end{align}
By virtue of the Sherman-Morrison-Woodbury formula, the matrices to be inverted are either diagonal or small ($p_k \times p_k$). The right hand side becomes:
\begin{align}
\bm{r}_1 - \sum_{k=2}^K \bm{B}_k \left( \bm{A}_{k1}\bm{W}_1\bm{A}_{k1}^T + \bm{E}_k \bm{W}_k \bm{E}_k^T \right)^{-1} \bm{r}_k, \label{eq:reduced_rhs}
\end{align}
where $\bm{r}_k$ is the part of the right hand side of \eqref{eq:norm} that corresponds to block $k$.
Both computing the right hand side and recovering $\Delta\bm{y}_k$ from $\Delta\bm{y}_1$ can also be implemented via the Sherman-Morrison-Woodbury formula.

We continue with an analysis of the number of flops to compute the coefficient matrix, shown in Table \ref{tbl:flops}. It is assumed that $p_k \leq m_k$. The full system \eqref{eq:norm:coeff} can only be computed in one way. The reduced system \eqref{coeff:reduced} can be computed in multiple ways. We show the number of flops for an efficient method, e.g., it avoids computing and storing the full matrix $\bm{B}_k \in \mathbb{R}^{m_1 \times m_k}$. To compute the second term in \eqref{coeff:reduced}, we use the original formula for $\bm{B}_k$ and expand the term to:
\begin{align}
&\bm{A}_{11}\bm{W}_1\bm{A}_{k1}^T  \left(  \bm{E}_k \bm{W}_k \bm{E}_k^T \right)^{-1} \bm{A}_{k1} \bm{W}_1 \bm{A}_{11}^T +
\bm{A}_{11}\bm{W}_1\bm{A}_{k1}^T  \left(  \bm{E}_k \bm{W}_k \bm{E}_k^T \right)^{-1}\bm{E}_k \bm{W}_k \bm{A}_{1k}^T + \notag \\
&\qquad \qquad \bm{A}_{1k}\bm{W}_k\bm{E}_k^T  \left(  \bm{E}_k \bm{W}_k \bm{E}_k^T \right)^{-1} \bm{A}_{k1} \bm{W}_1 \bm{A}_{11}^T +
\bm{A}_{1k}\bm{W}_k\bm{E}_k^T  \left(  \bm{E}_k \bm{W}_k \bm{E}_k^T \right)^{-1} \bm{E}_k \bm{W}_k \bm{A}_{1k}^T. \label{eq:bb}
\end{align}
As $\bm{A}_{k1}$ is assumed to have a limited number of nonzero columns, the first term in \eqref{eq:bb} can be computed as $\sum_j \left((\bm{A}_{k1})_j^T  \left(  \bm{E}_k \bm{W}_k \bm{E}_k^T \right)^{-1} (\bm{A}_{k1})_j (\bm{W}_1)_j^2\right) (\bm{A}_{11})_j (\bm{A}_{11})_j^T$, where $j$ runs over the nonzero columns of $\bm{A}_{k1}$. This is a scalar multiplied with a rank one matrix. The second and third term differ merely by a transposition. They are $\bm{O}$ if the columns of $\bm{A}_{11}$ that correspond to the nonzero columns of $\bm{A}_{k1}$ are zero. Otherwise their computational cost is similar to that of $\bm{A}_{11}\bm{A}_{k1}^T\bm{E}_k \bm{E}_k^T \bm{A}_{1k}^T$, which can be computed as $\bm{A}_{11}((\bm{A}_{k1}^T(\bm{E}_k \bm{E}_k^T) )\bm{A}_{1k}^T)$ by propagating the nonzero rows of $\bm{A}_{k1}^T$. The computational cost for the last term in \eqref{eq:bb} depends on the nonzero structure of $\bm{E}_k^T \bm{E}_{k}$. The number of flops for this step in Table \ref{tbl:flops} assumes that $\bm{E}_k^T \bm{E}_{k}$ is dense. When that product is sparse, which will be the case in many applications, the last term in \eqref{eq:bb} can be computed more efficiently than reported in the table. The computation of the final terms in \eqref{coeff:reduced} is divided into four steps, where the matrix $\bm{B}_k$ is expanded again.

There may be a discrepancy between $p_k$ and the rank of $\bm{A}_{k1}$. In that case there exists a $\bm{V}_k$ with fewer columns than $p_k$ such that $\bm{A}_{k1}\bm{W}_1\bm{A}_{k1}^T = \bm{V}_k\bm{V}_k^T$, and the number of flops reported for $\bm{C}_k$ and $\bm{D}_k$ in Table \ref{tbl:flops} can be reduced at the expense of factorizing $\bm{A}_{k1}\bm{W}_1\bm{A}_{k1}^T$.

\begin{table}
\caption{Number of flops required to compute the coefficient matrix, assuming $p_k \leq m_k$.}\label{tbl:flops}
\hspace{-2cm}
\begin{tabular}{llll}
\toprule
\multicolumn{2}{c}{Full system \eqref{eq:norm:coeff}} & \multicolumn{2}{c}{Reduced system \eqref{coeff:reduced}} \\
\cmidrule(lr){1-2} \cmidrule(lr){3-4}
Formula & Flops & Formula & Flops \\
\midrule
$\sum_{k=1}^K \bm{A}_{1k} \bm{W}_k \bm{A}_{1k}^T$ & $\mathcal{O}(m_1^2 n)$
& $\sum_{k=1}^K \bm{A}_{1k} \bm{W}_k \bm{A}_{1k}^T$ & $\mathcal{O}(m_1^2 n)$ \\
$\bm{A}_{11}\bm{W}_1\bm{A}_{k1}^T+\bm{A}_{1k}\bm{W}_k\bm{E}_k^T$ & $\mathcal{O}(m_1 m_k n_1)$ 
& $\bm{A}_{11}\bm{W}_1\bm{A}_{k1}^T  \left(  \bm{E}_k \bm{W}_k \bm{E}_k^T \right)^{-1} \bm{A}_{k1} \bm{W}_1 \bm{A}_{11}^T$ & $\mathcal{O}(m_1^2m_kp_k)$ \\
$\bm{A}_{k1}\bm{W}_1\bm{A}_{k1}^T + \bm{E}_k \bm{W}_K \bm{E}_k^T$ & $\mathcal{O}(m_k^2)$ 
& $\bm{A}_{11}\bm{W}_1\bm{A}_{k1}^T  \left(  \bm{E}_k \bm{W}_k \bm{E}_k^T \right)^{-1}\bm{E}_k \bm{W}_k \bm{A}_{1k}^T$ & $\mathcal{O}(m_1 m_k p_k)$ \\
& & $\bm{A}_{1k}\bm{W}_k\bm{E}_k^T  \left(  \bm{E}_k \bm{W}_k \bm{E}_k^T \right)^{-1} \bm{E}_k \bm{W}_k \bm{A}_{1k}^T$ & $\mathcal{O}(m_k n_k^2\!+\!m_1^2 m_k)$  \\
& & $\bm{C}_k=\left( \bm{I}_{p_k} + \bm{V}_k^T \left(  \bm{E}_k \bm{W}_k \bm{E}_k^T  \right)^{-1} \bm{V}_k \right)^{-0.5}$ & $\mathcal{O}(m_k^2p_k)$ \\
& & $\bm{D}_k=\left(  \bm{E}_k \bm{W}_k \bm{E}_k^T \right)^{-1} \bm{V}_k^T \bm{C}_k$ & $\mathcal{O}(m_k^2p_k)$ \\
& & $\bm{G}_k = \bm{A}_{11}\bm{W}_1\bm{A}_{k1}^T \bm{D}_k+\bm{A}_{1k}\bm{W}_k\bm{E}_k^T \bm{D}_k$ & $\mathcal{O}((m_1\!+\!m_k)(n_1\!+\!n_k)p_k)$ \\
& & $\bm{G}_k \bm{G}_k^T$ & $\mathcal{O}(m_1^2 p_k)$ \\
\bottomrule
\end{tabular}
\end{table}

The time to compute the coefficient matrix \eqref{coeff:reduced} can be related to the structure of the CPL constraint \eqref{eq:piecewise}. Let $j$ denote the nonzero column in $\bm{A}_{k1}$. The first term in \eqref{eq:bb} simplifies to
$$\left( (\bm{A}_{k1})_j^T \left(  \bm{E}_k \bm{W}_k \bm{E}_k^T \right)^{-1} (\bm{A}_{k1})_j \right) (\bm{A}_{11})_j (\bm{A}_{11})_j^T,$$
which is a scalar multiplied with a rank one matrix. The second and third term in \eqref{eq:bb} are $\bm{O}$. In the fourth term, $\bm{E}_k^T \bm{E}_{k}$ is block diagonal with size $L$. The smaller the block size, the easier it is to compute the final term in \eqref{eq:bb}, the easiest cases being the maximum of 0 and a linear function, or the maximum of a function and the negative of that function (to model the absolute value).

The cost of solving \eqref{eq:norm} directly is $(1/3)m^3 + 2m^2$ flops, while the cost of solving the reduced system is $(1/3)m_1^3 + 2m_1^2$ flops. The savings come at the expense of row reduction steps for the coefficient matrix and the right hand side, which is a minor added cost in big-$\mathcal{O}$ sense as long as $p_k$ is limited.

\subsection{Comparison with other IPMs}\label{sec:comparison}
\citet{castro2016interior} developed an interior point method for a structure similar to \eqref{coeff:structure}, where he also eliminated all but one block. The key difference is that he did not consider a specific structure for the blocks on the diagonal. Therefore, the formula to substitute out $\Delta \bm{y}_k$ \eqref{eq:deltay} was not simplified with the Sherman-Morrison-Woodbury formula, but substituted straight into the reduced coefficient matrix. The reduced system is solved with a preconditioned conjugate gradient method. The preconditioner is created with a power series which converges slowly for linear optimization problems as \citet{castro2016interior} notes. The Sherman-Morrison-Woodbury formula avoids these numerical difficulties.

Instead of performing row reduction steps on \eqref{coeff:structure}, it is also possible to factorize the full matrix \eqref{coeff:structure}. Typically the rows and columns are first permuted in such a way that its Cholesky factorization becomes sparse \citep[\S 9.7.2]{Boydcvx}. For this particular matrix, such a permutation can yield a block arrow matrix. The Cholesky factor then has a bordered form \citep[see, e.g.,][]{gondzio2009exploiting} plus a low rank corrector:
\begin{align*}
\begin{pmatrix}
\bm{D}_2 + \bm{R}_2 & \hdots   & \bm{O}    & (\bm{AA}^T)_{12}^T \\
\vdots              & \ddots   & & \\
\bm{O}              &          & \bm{D}_K + \bm{R}_K & (\bm{AA}^T)_{1K}^T \\
(\bm{AA}^T)_{12}    &          & (\bm{AA}^T)_{1K} & (\bm{AA}^T)_{11}
\end{pmatrix}
= \bm{L}\bm{L}^T + \bm{R}\bm{R}^T,
\end{align*}
\begin{align*}
\text{with } \bm{L} = \begin{pmatrix}
\bm{L}_2 \\
& \ddots \\
& & \bm{L}_K \\
\bm{L}_{K2} & \hdots & \bm{L}_{KK} & \bm{L}_c
\end{pmatrix}
\text{, }
\bm{L}_k = \bm{D}_k^\frac{1}{2}
\text{, }
\bm{L}_{Kk} = \bm{D}_k^{-\frac{1}{2}} (\bm{AA}^T)_{1k}
\text{, }
\end{align*}
and $\bm{L}_c$ is the Cholesky factor of $(\bm{AA}^T)_{11} - \sum_{k=2}^K (\bm{AA}^T)_{1k} \bm{D}_i^{-1} (\bm{AA}^T)_{1k}$. Applying a low rank update to the Cholesky factor creates burdensome fill-in, so this method should employ the Sherman-Morrison-Woodbury formula as well, making it mathematically equivalent to our approach. Although \eqref{coeff:structure} has never been factorized this way to the best of our knowledge, it combines the ideas in Section 4.1 and 4.2 of \cite{gondzio2003parallel}.

\section{Numerical examples}\label{sec:examples}
\subsection{Netlib}
To get a sense of the prevalence of the exploitable structure and the efficacy of Alg.~\ref{alg:structure}, we analyze the Netlib/lp test set \citep{gay1985electronic}, which has been made publicly available in Matlab format by Tim Davis (Texas A\&M University, College Station, TX, USA). The set consists of 146 problems that were mostly collected during the late 80s and early 90s.

Alg.~\ref{alg:structure} detects structure in 90\% of the problems (in 64\% of the primal problems and 87\% of the dual problems). For 20\% of these problems, the dimension of the normal equations can be reduced by at least 50\% (in 8\% of the primal problems and in 13\% of the dual problems).

Further analysis reveals that $\bm{A}_{k1}$ often has rank 0, which is not as interesting because then the coefficient matrix \eqref{eq:norm:coeff} is sparse already. We therefore slightly modified the algorithm to detect only the cases where $\bm{A}_{k1}$ is nonzero. The structure is still detected in 36\% of the problems (in 4\% of the primal problems and 34\% of the dual problems). The problems where the structure is most prominent are {\sc lpi\_cplex1}, {\sc scrs8} and {\sc modszk1}, where the dimension of the normal equations can be reduced by 50\%, 38\% and 36\%, respectively.

\subsection{Radiation therapy} \label{sec:radiation}
To quantitatively analyze the performance of the proposed method, we use the optimization model and data from the proton cases in the publicly available data set TROTS \citep{breedveld2017data}. They are all head and neck cases which share the same set of relevant structures, objectives and constraints, but differ in the location of the tumor, the pencil beam placement and objective weights. The sizes of the cases are provided in Table \ref{tbl:trots}. The number of pencil beams varies from 990 to 2761 which is representative for proton as well as photon cases, although larger tumors or a smaller spot size lead to an increase in the number of pencil beams. The dimension of the normal equations is 110,047 on average, which is reduced by our method to 1,829 (an average reduction of 98.4\%).

The objectives and constraints are all on the minimum, maximum and mean dose. Some objectives and constraints are {\it robust} against nine scenarios (a nominal, undershoot and overshoot scenario, and set-up errors in the $x$/$y$/$z$ direction) in the sense that constraints should hold for each scenario and objectives are optimized against the worst case \citep{unkelbach2018robust}. Mathematically, the optimization problems for the proton cases of TROTS can be stated as:
\begin{align*}
\min_{\bm{x},\bm{z}}  \quad & \sum_{s \in O} w_s z_s \\
\text{s.t.} \quad & \sum_j D_{ij} x_j \leq z_s \quad \forall i \in I(s), s \in O \\
                  & l_s \leq \sum_j D_{ij} x_j \leq u_s \quad \forall i \in I(s), s \in C \\
                  & \bm{x} \geq \bm{0},
\end{align*}
where $C$ is the set of constrained structures, $O$ is the set of structures for which the maximum dose is optimized with weight $w_o$, $I(s)$ is the set of voxels in structure $s$, and $l_s$ and $u_s$ are lower- and upper bounds on the dose in structure $s$. Note that a constraint on mean dose fits this framework by adding a single voxel structure, and a robust constraint simply increases the number of voxels for that constraint by a factor of 9. Deviating from the original problem specification, we modify the constraints on $s\in C$ to allow a slight constraint violation, as long as the average violation does not exceed 0.1 Gy:
\begin{align*}
\frac{1}{|I(s)|}\sum_{i \in I(s)} \max\left\{0, \sum_j D_{ij} x_j - u_s\right\} \leq 0.1 \quad \forall s \in C \\
\frac{1}{|I(s)|}\sum_{i \in I(s)} \max\left\{0, l_s - \sum_j D_{ij} x_j\right\} \leq 0.1 \quad \forall s \in C.
\end{align*}
This modification introduces exploitable structure, and reduces the impact of a small subset of voxels on the overall plan. Moreover, a single voxel cannot render the problem infeasible, thus, it is easier to formulate a problem that is feasible compared to using minimum or maximum dose constraints, which is useful to narrow down the search space for Pareto navigation. Mean overdose
and underdose constraints are by no means exotic: they are considered extremely useful for treatment planning \citep{craft2008many}. The resulting formulation for case 1 is shown in Table \ref{tbl:trots1}.

We aimed for a complete and equitable comparison against leading algorithms, and investigated the differences in objective value, solution time and plan quality. In Section \ref{eq:comparison_cplex} we compare implementations of our algorithm with the state-of-the-art interior point solver CPLEX. In Section \ref{sec:dosimetric_comparison} we compare the plan quality of our exact optimization algorithm with L-BFGS, which is used in many treatment planning systems, to verify if a deterioration in objective value has an effect on plan quality. The following algorithms were not included in the comparison, for reasons outlined below.

The alternating direction method of multipliers (ADMM) is an algorithm that splits up the optimization problem into a subproblem in the dose domain, a subproblem in the fluence domain, and a coupling problem, and was successfully applied to radiation therapy by \citet{ungun2019real}. We have tried various implementations of ADMM (i.e., POGS, OSQP, Optkit) but have not been able to obtain a reasonable solution for the TROTS proton set, even with a day of computation time. For example, OSQP was unable to find a solution where $x$ did not have several highly negative components. For Optkit we modified TROTS case 4, which is the smallest case, to $\min_{\bm{x} \geq \bm{0}} \{ ||\bm{Dx}-\bm{p}||_1 \}$, where $\bm{p}=\bm{Dx^*}$ and $\bm{x^*}$ is an optimal solution to the original TROTS problem. Even after downsampling the case to 100 voxels and running Optkit for a million iterations, the objective value is larger than $10^3$ while the optimum is $0$ by construction, which means the average deviation per voxel is 10 Gy. For comparison, \citet{ungun2019real} used Optkit to optimize the 1-norm and reported near optimal solutions for a case with 268,228 voxels and 330 beamlets after 3117 iterations. We could not resolve the discrepancy in results because they do not share their data. ADMM has a penalty parameter $\rho$ that is adjusted automatically by Optkit. We have performed additional experiments where ADMM was run for a million iterations per fixed value for $\rho$, which was varied between $10^{-3}$ and $1.04 \cdot 10^3$ in 20\% increments, but no value could improve on the automatic adjustment procedure.

ART3+O is an optimization algorithm built around the feasibility seeking method ART3+. An example of a feasibility problem is ``is an objective value of 30 Gy attainable?'', which can be answered by finding a point in a polyhedral set (or concluding that the set is empty). ART3+ iterates over the inequalities that define the polyhedron, and for each inequality that is violated, it reflects the current solution into the inequality, until all constraints are satisfied. It is proven to find a feasible point in a finite number of steps if such a point exists, and was successfully applied to proton therapy \citep{chen2010fast}. To avoid thresholding issues for concluding that a polyhedron is empty, we called ART3+ to find a feasible solution with an objective value of 50 Gy, and reduced that target objective by 0.1 Gy at a time, reusing the previous solution as the starting point. After $10^{12}$ iterations, the best solution had an objective value of 35.3 Gy, and it took approximately $55\cdot 10^9$ iterations to get there from the solution with a value of 35.4 Gy. Since it took around a day to obtain this result and the objective value is still considerably worse than what we obtained with L-BFGS, we did not include this algorithm in our comparison.

\begin{table}
\caption{Description of the TROTS data set and results from L-BFGS (IPOPT) after 500 iterations. Infeasibility is summed over the constraints.}\label{tbl:trots}
\begin{tabular}{lrrrrcc}
\toprule
Case & \multicolumn{2}{c}{Instance size} & \multicolumn{2}{c}{Dimension} & Suboptimality & Infeasibility \\
     \cmidrule(lr){2-3} \cmidrule(lr){4-5}
   & Beamlets & Voxels & $\bm{AA}^T$ & $(\bm{AA}^T)_{11}$ & (Gy) & (Gy) \\
\midrule
 1 & 1080 & 332,704 &  88,744 & 1107 & 1.84 & 0.00 \\
 2 & 2062 & 350,495 & 114,627 & 2089 & 5.38 & 0.01 \\
 3 & 2016 & 358,314 & 121,763 & 2036 & 4.69 & 0.00 \\
 4 &  990 & 331,009 &  81,694 & 1017 & 1.99 & 0.00 \\
 5 & 1280 & 328,859 &  83,668 & 1307 & 2.76 & 0.00 \\
 6 & 2618 & 387,416 & 138,406 & 2645 & 5.92 & 0.00 \\
 7 & 1881 & 335,742 & 101,052 & 1908 & 3.23 & 0.00 \\
 8 & 1109 & 331,747 &  88,950 & 1136 & 1.62 & 0.00 \\
 9 & 2761 & 382,076 & 141,102 & 2788 & 8.18 & 0.02 \\
10 & 2406 & 356,893 & 119,740 & 2433 & 6.73 & 0.01 \\
11 & 1803 & 391,067 & 150,188 & 1830 & 4.57 & 0.09 \\
12 & 1344 & 332,162 &  98,240 & 1371 & 3.41 & 0.00 \\
13 & 2252 & 359,409 & 118,637 & 2279 & 3.76 & 0.05 \\
14 & 2483 & 373,316 & 132,163 & 2510 & 8.87 & 0.03 \\
15 & 1962 & 333,145 & 102,793 & 1989 & 4.20 & 0.00 \\
16 & 2266 & 346,511 &  97,428 & 2293 & 6.42 & 0.05 \\
17 & 1400 & 339,482 & 104,140 & 1427 & 2.73 & 0.00 \\
18 & 1587 & 363,199 & 117,487 & 1614 & 3.08 & 0.00 \\
19 & 1675 & 359,766 & 117,157 & 1702 & 6.45 & 0.46 \\
20 & 1079 & 315,253 &  82,974 & 1106 & 2.15 & 0.00 \\
\midrule
mean & 1803 & 350,428 & 110,048 & 1829 & 4.40 & 0.04 \\
\bottomrule
\end{tabular}
\end{table}

\subsubsection{Algorithmic comparison}\label{eq:comparison_cplex}
We compare three different IPM implementations. The first is CPLEX 12.8, which represents the state-of-the-art in general purpose linear optimization algorithms. The second is an IPM written in Python from scratch for the purpose of this paper that detects and exploits structure. The third is Nymph, a proprietary IPM developed by the author specifically for radiation therapy that exploits structure in the same way as the Python version. The Python IPM has disadvantages that can be avoided in a lower level language: (1) all operations are single-threaded and do not take advantage of a multi-core cpu, (2) sparse matrix multiplications perform symbolic multiplication to determine the sparsity structure of the product in each iteration, even though the sparsity structure does not change between iterations or is dense, (3) each matrix element of a product is computed explicitly, even for a symmetric matrix, (4) intermediate calculation results are stored separately even if they can be added directly to an existing matrix, and (5) memory management is inefficient due to garbage collection. The Python program also does not have higher order correctors. CPLEX was run in barrier mode (without crossover) on the dual problem, which was found to be the fastest setting. The stopping criteria of Nymph were set to a duality gap of $10^{-8}$ Gy and primal and dual infeasibility of $||Ax-b|| \leq 10^{-5}(1+||b||)$ and $||A^Ty+s-c|| \leq 10^{-5}(1+||c||)$. The most stringent criterion is the one for the duality gap, and a setting of $10^{-8}$ matches with CPLEX. The algorithms were run on a dual socket Intel Xeon E5-2687W v3 with a total of 20 physical cpu cores.

The first comparison is based on a single IPM iteration with the predictor-corrector method disabled. We timed the actual time per iteration, ignoring preprocessing. This is the most accurate way to compare the gains from exploiting structure, since algorithmic differences are mostly eliminated. The second comparison is based on total runtime. Finally, we consider the number of iterations.

The results in Table \ref{tbl:results} show that Python is an order of magnitude slower than CPLEX, despite exploiting structure. This is clear both in the time per iteration and in the total runtime. Python and CPLEX need approximately the same number of iterations. Apparently there is little benefit from higher order corrections in CPLEX, which was also concluded by \citep{Breedveld2016} for the unmodified TROTS set. CPLEX is an order of magnitude faster than Python, which is not unexpected given the factors outlined above. In turn, Nymph is an order of magnitude faster than CPLEX. Compared to the original TROTS formulation, Nymph slows down less than 45\% due to the mean overdose constraints whereas CPLEX slows down by a factor of 3.

\begin{table}
\caption{Runtimes and statistics of CPLEX and two implementations of a structure exploiting interior point method.}\label{tbl:results}
\begin{tabular}{lllllllll}
\toprule
Case & \multicolumn{3}{c}{Time per iteration (s)} & \multicolumn{3}{c}{Total time (s)} & \multicolumn{2}{c}{Number of iterations} \\
\cmidrule(lr){2-4} \cmidrule(lr){5-7} \cmidrule(lr){8-9}
& CPLEX & Python & Nymph & CPLEX & Python & Nymph & CPLEX & Python \\
\midrule
1 & 2.2 & 14.5 & 0.2 & 316 & 2387 & 26 & 109 & 130 \\
2 & 4.1 & 30.8 & 0.4 & 636 & 4893 & 54 & 154 & 133 \\
3 & 3.6 & 27.9 & 0.4 & 831 & 5569 & 66 & 206 & 162 \\
4 & 2.1 & 13.1 & 0.2 & 275 & 2359 & 32 & 106 & 140 \\
5 & 3.5 & 18.4 & 0.3 & 416 & 3919 & 40 & 133 & 170 \\
6 & 4.1 & 33.5 & 0.5 & 703 & 5710 & 65 & 142 & 135 \\
7 & 2.6 & 23.4 & 0.3 & 453 & 3933 & 41 & 150 & 138 \\
8 & 2.1 & 14.1 & 0.2 & 250 & 2137 & 28 & 107 & 119 \\
9 & 4.5 & 35.5 & 0.5 & 744 & 5784 & 60 & 149 & 129 \\
10 & 4.3 & 34.0 & 0.4 & 666 & 5703 & 59 & 138 & 137 \\
11 & 4.0 & 28.0 & 0.3 & 477 & 4127 & 40 & 106 & 124 \\
12 & 3.2 & 21.0 & 0.3 & 520 & 2912 & 33 & 102 & 112 \\
13 & 3.8 & 33.2 & 0.4 & 630 & 6046 & 55 & 151 & 148 \\
14 & 4.5 & 34.4 & 0.4 & 661 & 5024 & 47 & 113 & 118 \\
15 & 3.1 & 24.4 & 0.3 & 347 & 3165 & 34 & 100 & 105 \\
16 & 4.1 & 29.8 & 0.4 & 623 & 4345 & 42 & 136 & 117 \\
17 & 2.3 & 17.4 & 0.2 & 301 & 2259 & 26 & 107 & 103 \\
18 & 3.0 & 21.3 & 0.3 & 341 & 2923 & 31 & 99 & 109 \\
19 & 3.1 & 24.5 & 0.3 & 405 & 4657 & 43 & 114 & 152 \\
20 & 3.6 & 17.3 & 0.2 & 315 & 2650 & 29 & 100 & 118 \\
\midrule
mean & 3.4 & 24.8 & 0.3 & 495.6 & 4025.1 & 42.5 & 126.1 & 130.0 \\
\bottomrule
\end{tabular}
\end{table}

\subsubsection{Dosimetric comparison}\label{sec:dosimetric_comparison}
We compare Nymph to the L-BFGS solver IPOPT \citep{wachter2006implementation}. We ran IPOPT for 500 iterations, which is an order of magnitude more than what is typically done clinically. The goal of this comparison is to show the dosimetric differences between L-BFGS and an exact optimization algorithm. Therefore, timing information is not shown. It should be noted that highly optimized implementations of L-BFGS can run in seconds \citep{ziegenhein2013performance}. Compared to the previous section, the termination criterion of Nymph was relaxed to allow a duality gap of 0.1 Gy, which is a trade-off between solution time and quality.

We use the log-sum-exp function with scaling parameter $\varepsilon = 10^{-3}$ as a conservative approximation to the maximum function \citep{fredriksson2014critical}:
\begin{align*}
\max_{i \in I(s)} \left\{\sum_j D_{ij} x_j \right\} & \approx \varepsilon \log\left( \sum_{i \in I(s)} e^{\sum_j D_{ij} x_j / \varepsilon}\right) \text{, and} \\
\sum_{i \in I(s)} \max\left\{0, \sum_j D_{ij} x_j - u_s\right\} & \approx \sum_{i \in I(s)} \varepsilon \log\left(1 + e^{\left(\sum_j D_{ij} x_j-u_s\right) / \varepsilon}\right).
\end{align*}
Therefore, instead of having one linear constraint per voxel resulting in hundreds of thousands of constraints, the number of constraints is now limited to 16. This step is necessary for L-BFGS to have good performance, although there is some variety in the choice of a nonlinear approximation.

\begin{figure}\label{fig:trots1}
\centering
\includegraphics{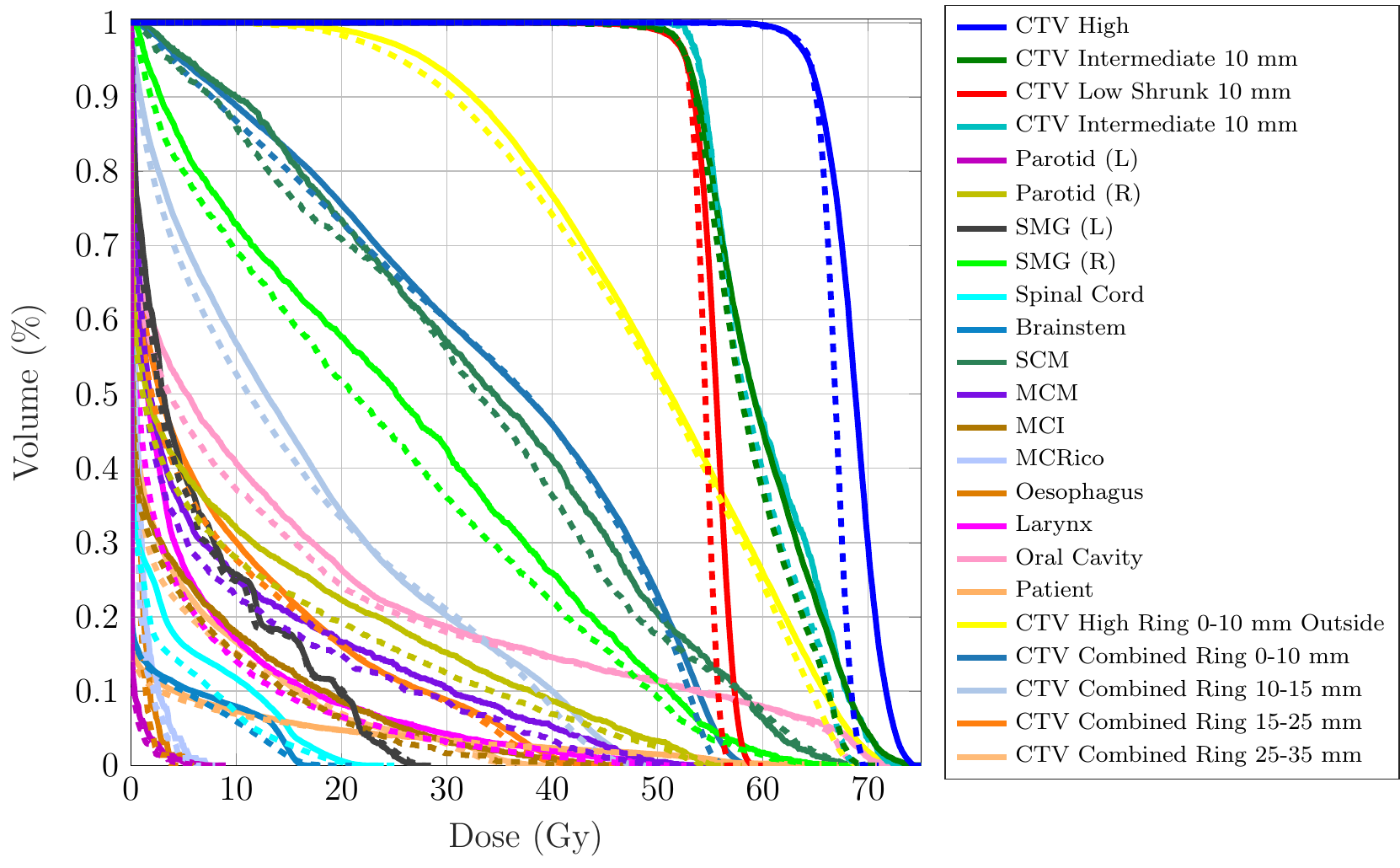}
\caption{DVH for TROTS case 1. Solid curves: L-BFGS. Dashed curves: Nymph.}
\end{figure}

The solutions of IPOPT and Nymph are compared in Table \ref{tbl:trots1} and Figure \ref{fig:trots1} for the first case, and summarized in Table \ref{tbl:trots} for all twenty cases. Table \ref{tbl:trots1} shows that both plans satisfy the constraints, but that Nymph has an objective value that is 1.84 Gy smaller. Since the weights of the maximum and mean dose objectives sum to 1, the (weighted) average improvement per objective value is 1.84 Gy. Most of the improvements are in the maximum dose in the CTV and the mean dose to many avoidance structures. Indeed the DVH in Figure \ref{fig:trots1} shows that the solution of Nymph has better target homogeneity and that the mean dose is improved for many structures. Across all 20 cases, the improvement in objective value is 4.40 Gy on average and ranges between 1.62 Gy (case 8) and 8.87 Gy (case 14).

The IPOPT solutions only have a few violated constraints. The only case that is severely affected is case 19 where two constraints are violated by more than 100\%, which results in a clear deterioration of the corresponding DVHs. Constraint violations for this case could be avoided by increasing the number of iterations to 1400. Nymph does not have a single constraint violation because that is part of the termination criterion.

To investigate the effects of the nonlinear reformulation, we tested the following model variations on the first three TROTS cases:
\begin{enumerate}
\item Replacing the constrained problem with an unconstrained problem using the sixteen optimal Lagrange multipliers of the constraints. Eliminating the constraints takes away a key difference between implementations of L-BFGS, e.g., augmented Lagrangian and interior point methods. The objective function is a weighted sum of maximum dose, mean overdose and mean underdose functions, which are still modeled via log-sum-exp. The objective value after 500 iterations was worse than what IPOPT achieved for the constrained problem (34.09, 51.67 and 53.60 for cases 1, 2 and 3, respectively).
\item Replacing log-sum-exp in the objective function with the $p$-norm. This resulted in plans with similar objective values as the log-sum-exp plans (32.71 and 32.40, 44.45 and 42.92, and 46.02 and 45.12 for $p=10$ and $p=20$, for cases 1, 2 and 3, respectively).
\item Modeling the maximum dose objectives as mean overdose objectives. The dose threshold of each objective was set to the maximum dose in an optimal solution. This improved the objective value (to 31.78, 43.12, 43.71 for cases 1, 2 and 3, respectively), but the plans are still suboptimal and the selected threshold level is not available in practice.
\item Comparing the algorithms on the original TROTS formulation, i.e., with maximum and minimum dose constraints. We tried three different convex reformulations of the constraints: (a) via log-sum-exp, (b) the mean overdose above or mean underdose below the maximum or minimum dose bound should be 0, and (c) same as method (b) except using mean squared overdose or underdose. None of the formulations was able to produce a decent treatment plan. The optimal objective values for the first three cases are 32.65, 42.92 and 44.09 Gy. Method (a) produced plans with an objective value of 36.98, 52.50 and 58.31 Gy and constraint violations of up to 2.5 Gy for the minimum CTV dose and maximum dose in the patient. Method (b) yielded plans with objective values of 35.49, 54.13 and 53.07 Gy, where the second plan was severely infeasible with violations of up to 4 Gy. Method (c) resulted in objective values of 46.14, 67.42 and 72.99 Gy, and all plans had a constraint violation of at least 4 Gy. In conclusion, all plans are highly suboptimal and most plans do not satisfy the constraints.
\end{enumerate}
From the observations above, it is clear that the problem formulation has an effect on the performance of L-BFGS. None of the tested formulations could achieve a treatment plan quality similar to Nymph.

We have additionally tried the two L-BFGS based solvers SNOPT \citep{gill2005snopt} and L-BFGS-B \citep{zhu1997algorithm}. Since L-BFGS-B only supports box constraints, we tested that algorithm with model variation 1 (using optimal Lagrange multipliers to put constraints in the objective). Each algorithm was run for 500 and 5000 iterations with limited BFGS memory, and for 500 iterations with full memory. The results in Table \ref{tbl:lbfgs} show that IPOPT and L-BFGS-B attain similar objective values after 500 and 5000 iterations in limited memory mode, which are better than SNOPT. Using the full memory improves the objective value at the expense of constraint violations. Even with 5000 iterations, all algorithms are at least 1 Gy away from optimality on average.

\begin{table}
\centering
\caption{Objective value for different implementations of L-BFGS. The column Mem indicates if limited or full memory was used for BFGS.}\label{tbl:lbfgs}
\begin{tabular}{lcrccclll}
\toprule
       &         &            & \multicolumn{3}{c}{Objective value}& \multicolumn{3}{c}{Total infeasibility} \\
       \cmidrule(lr){4-6} \cmidrule(lr){7-9}
Method & Mem & Iter. & Trots 1 & Trots 2 & Trots 3 & Trots 1 & Trots 2 & Trots 3 \\
\midrule
Nymph    &   &      & 30.3 & 40.7 & 41.6 \\
\addlinespace
IPOPT    & L &  500 & 32.2 & 46.1 & 46.3 & 0 & 0.001 & 0 \\
IPOPT    & F &  500 & 30.8 & 43.5 & 47.3 & 0.001 & 0.213 & 0.029 \\
IPOPT    & L & 5000 & 30.8 & 41.8 & 43.6 & 0 & 0 & 0 \\
\addlinespace
SNOPT    & L &  500 & 35.5 & 52.0 & 49.4 & 0 & 0 & 0 \\
SNOPT    & F &  500 & 30.2 & 49.8 & 49.5 & 0.232 & 0 & 0 \\
SNOPT    & L & 5000 & 32.6 & 46.9 & 45.6 & 0 & 0 & 0 \\
\addlinespace
L-BFGS-B & L &  500 & 32.6 & 45.5 & 47.1 \\
L-BFGS-B & F &  500 & 32.4 & 44.7 & 45.1 \\
L-BFGS-B & L & 5000 & 31.2 & 41.8 & 42.7 \\
\bottomrule
\end{tabular}
\end{table}

\section{Concluding remarks}
Virtually all treatment plans today are created with inexact optimization algorithms, i.e., without an optimality guarantee. It is widely believed that the radiation therapy optimization problem is relatively easy, and that many algorithms can produce treatment plans of high dosimetric quality. For intensity-modulated photon therapy (IMRT) this has been explained by analyzing the dose map \eqref{eq:dose} and the Hessian, which revealed that there are many (near) optimal solutions \citep{carlsson2006using,webb2003physical}. The objective function also does not always correlate well with the dosimetric quality, so being optimal does not necessarily translate into a better treatment \citep{alterovitz2006optimization,gorissen2013mixed}. However, there may be cases where the set of (near) optimal solutions is small, and objective functions can be formulated in a meaningful way to correlate with clinical endpoints. The proton cases from the TROTS data set seem to fit those criteria. We showed that most algorithms are unable to obtain an optimal solution, and that suboptimality is not merely a technicality, but that the dose distribution is severely affected. The issues are more apparent in the unmodified TROTS cases, where the feasible region is smaller than after modification, suggesting that the size of the set of (near) optimal solutions is indeed influential for L-BFSG. Optimization algorithms with an optimality guarantee therefore have a tangible benefit. Our work shows that such an algorithm can be implemented in a way that is fast enough for clinical use.

\begin{table}
\caption{Formulation and dose statistics for TROTS case 1}\label{tbl:trots1}
\hskip-1.0cm\begin{tabular}{llll}
\toprule
Description & Bound / weight & IPOPT & Nymph \\
\midrule
Robust mean underdose CTV High below 64.68 Gy & 0.1 Gy & 0.10 & 0.10\\
Mean underdose CTV Intermediate 10 mm below 52.92 Gy & 0.1 Gy & 0.07 & 0.09\\
Mean underdose CTV Low Shrunk 10 mm below 52.92 Gy & 0.1 Gy & 0.10 & 0.10\\
Mean overdose Parotid (L) above 69.96 Gy & 0.1 Gy & 0.00 & 0.00\\
Mean overdose Parotid (R) above 69.96 Gy & 0.1 Gy & 0.00 & 0.00\\
Mean overdose SMG (L) above 69.96 Gy & 0.1 Gy & 0.00 & 0.00\\
Mean overdose SMG (R) above 69.96 Gy & 0.1 Gy & 0.00 & 0.00\\
Mean overdose SCM above 69.96 Gy & 0.1 Gy & 0.00 & 0.00\\
Mean overdose MCM above 69.96 Gy & 0.1 Gy & 0.00 & 0.00\\
Mean overdose MCI above 69.96 Gy & 0.1 Gy & 0.00 & 0.00\\
Mean overdose MCRico above 69.96 Gy & 0.1 Gy & 0.00 & 0.00\\
Mean overdose Oesophagus above 69.96 Gy & 0.1 Gy & 0.00 & 0.00\\
Mean overdose Larynx above 69.96 Gy & 0.1 Gy & 0.00 & 0.00\\
Mean overdose Oral Cavity above 69.96 Gy & 0.1 Gy & 0.01 & 0.00\\
Mean overdose CTV Intermediate 10 mm above 100 Gy & 0.1 Gy & 0.00 & 0.00\\
Mean overdose Patient above 69.96 Gy & 0.1 Gy & 0.00 & 0.00\\
Robust minimize maximum dose CTV High & $3.4\cdot 10^{-2}$ & 74.19 & 69.42\\
Minimize maximum dose CTV Intermediate 10 mm & $2.6\cdot 10^{-3}$ & 74.95 & 69.49\\
Minimize maximum dose CTV Low Shrunk 10 mm & $2.7\cdot 10^{-1}$ & 58.88 & 56.59\\
Minimize maximum dose CTV High Ring 0-10 mm Outside & $5.2\cdot 10^{-3}$ & 74.52 & 69.32\\
Minimize maximum dose CTV Combined Ring 0-10 mm & $8.5\cdot 10^{-2}$ & 58.68 & 55.93\\
Minimize maximum dose CTV Combined Ring 10-15 mm & $8.3\cdot 10^{-3}$ & 54.22 & 48.80\\
Robust minimize mean dose Parotid (L) & $1.2\cdot 10^{-9}$ & 0.36 & 0.28\\
Robust minimize mean dose Parotid (R) & $3.3\cdot 10^{-1}$ & 12.52 & 11.08\\
Robust minimize mean dose SMG (L) & $5.8\cdot 10^{-2}$ & 8.60 & 8.31\\
Robust minimize mean dose SMG (R) & $4.0\cdot 10^{-2}$ & 31.02 & 29.20\\
Minimize maximum dose Spinal Cord & $1.6\cdot 10^{-3}$ & 22.72 & 23.53\\
Minimize maximum dose Brainstem & $7.3\cdot 10^{-6}$ & 16.99 & 19.82\\
Robust minimize mean dose SCM & $1.3\cdot 10^{-2}$ & 38.62 & 36.66\\
Robust minimize mean dose MCM & $1.3\cdot 10^{-2}$ & 12.24 & 11.09\\
Robust minimize mean dose MCI & $3.0\cdot 10^{-2}$ & 8.61 & 7.37\\
Robust minimize mean dose MCRico & $1.8\cdot 10^{-2}$ & 2.50 & 2.02\\
Robust minimize mean dose Oesophagus & $4.0\cdot 10^{-2}$ & 2.53 & 1.85\\
Robust minimize mean dose Larynx & $3.1\cdot 10^{-2}$ & 8.53 & 7.28\\
Robust minimize mean dose Oral Cavity & $9.1\cdot 10^{-3}$ & 17.36 & 16.32\\
Minimize mean dose CTV High Ring 0-10 mm Outside & $2.7\cdot 10^{-4}$ & 50.66 & 49.68\\
Minimize mean dose CTV Combined Ring 0-10 mm & $2.0\cdot 10^{-3}$ & 34.35 & 33.59\\
Minimize mean dose CTV Combined Ring 10-15 mm & $7.6\cdot 10^{-4}$ & 16.22 & 15.38\\
Minimize maximum dose CTV Combined Ring 15-25 mm & $7.1\cdot 10^{-4}$ & 47.32 & 48.58\\
Minimize mean dose CTV Combined Ring 15-25 mm & $1.5\cdot 10^{-6}$ & 8.41 & 7.86\\
Minimize maximum dose CTV Combined Ring 25-35 mm & $7.1\cdot 10^{-4}$ & 40.69 & 43.20\\
Minimize mean dose CTV Combined Ring 25-35 mm & $6.4\cdot 10^{-3}$ & 4.33 & 3.95\\
Minimize monitor units & $7.3\cdot 10^{-7}$ & $3.6\cdot 10^{5}$ & $3.4\cdot 10^{5}$\\
Objective value & & 32.19 & 30.34\\
\bottomrule
\end{tabular}
\end{table}

\section*{Acknowledgment}
Supported in part by NIH U19 Grant 5U19CA021239-38. The author thanks J.~Gondzio (University of Edinburgh, Scotland) for discussions during the development of Nymph and anonymous referees for constructive feedback.

\appendix
\newpage
\section{Dosimetric comparison}
\begin{table}[H]
\caption{Formulation and dose statistics for TROTS case 1}
\hskip-1.0cm
\begin{tabular}{llll}
\toprule
Description & Bound / weight & IPOPT & Nymph \\
\midrule
Robust mean underdose CTV High below 64.68 Gy & 0.1 Gy & 0.10 & 0.10\\
Mean underdose CTV Intermediate 10 mm below 52.92 Gy & 0.1 Gy & 0.07 & 0.09\\
Mean underdose CTV Low Shrunk 10 mm below 52.92 Gy & 0.1 Gy & 0.10 & 0.10\\
Mean overdose Parotid (L) above 69.96 Gy & 0.1 Gy & 0.00 & 0.00\\
Mean overdose Parotid (R) above 69.96 Gy & 0.1 Gy & 0.00 & 0.00\\
Mean overdose SMG (L) above 69.96 Gy & 0.1 Gy & 0.00 & 0.00\\
Mean overdose SMG (R) above 69.96 Gy & 0.1 Gy & 0.00 & 0.00\\
\bottomrule
\end{tabular}
\end{table}
\begin{figure}
\centering
\includegraphics{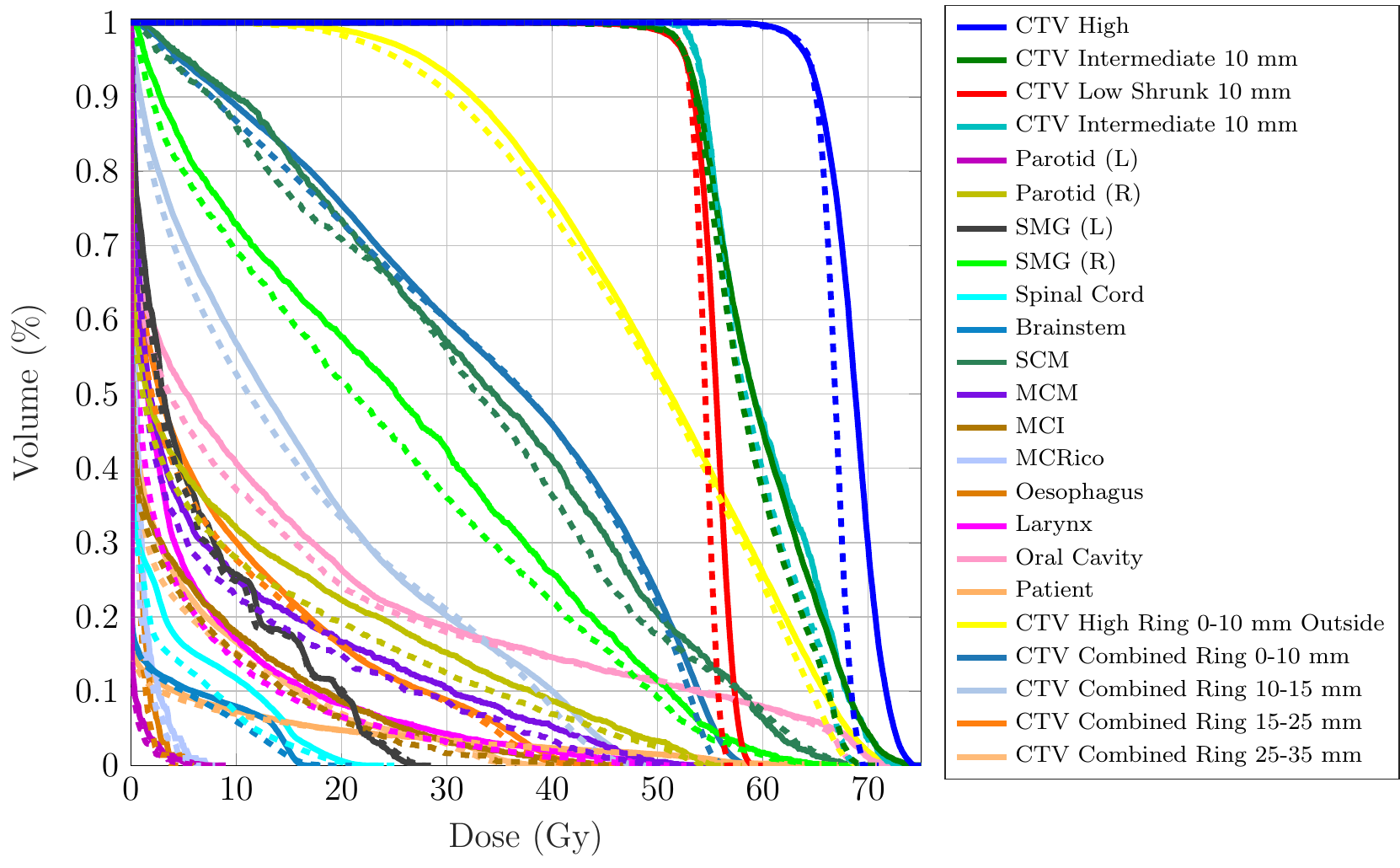}
\caption{DVH for TROTS case 1. Solid curves: L-BFGS. Dashed curves: Nymph.}
\end{figure}
\begin{figure}
\centering
\includegraphics{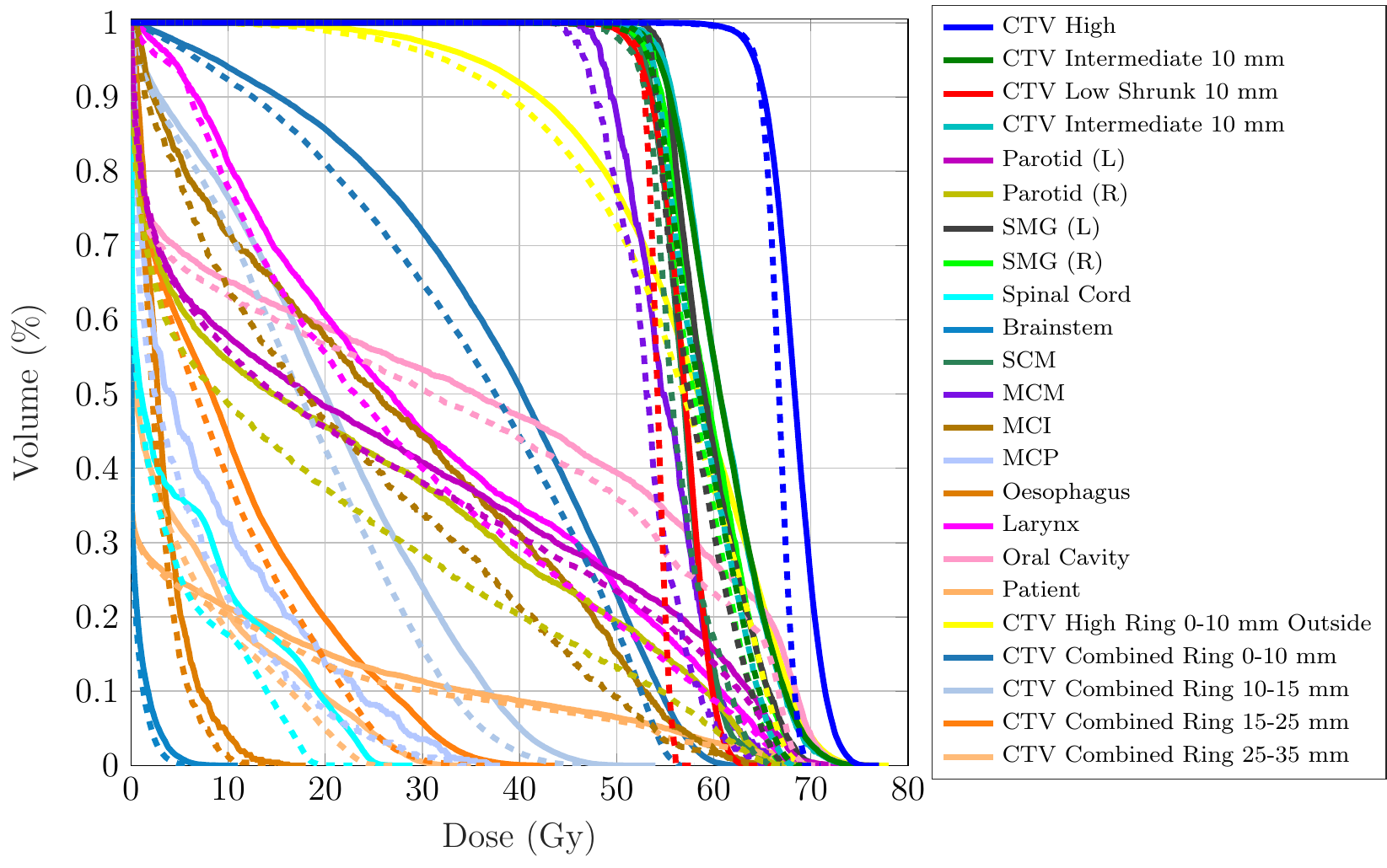}
\caption{DVH for TROTS case 2. Solid curves: L-BFGS. Dashed curves: Nymph.}
\end{figure}
\begin{table}
\caption{Formulation and dose statistics for TROTS case 2}
\hskip-1.0cm\begin{tabular}{llll}
\toprule
Description & Bound / weight & IPOPT & Nymph \\
\midrule
Robust mean underdose CTV High below 64.68 Gy & 0.1 Gy & 0.10 & 0.10\\
Mean underdose CTV Intermediate 10 mm below 52.92 Gy & 0.1 Gy & 0.03 & 0.06\\
Mean underdose CTV Low Shrunk 10 mm below 52.92 Gy & 0.1 Gy & 0.11 & 0.10\\
Mean overdose Parotid (L) above 69.96 Gy & 0.1 Gy & 0.01 & 0.00\\
Mean overdose Parotid (R) above 69.96 Gy & 0.1 Gy & 0.00 & 0.00\\
Mean overdose SMG (L) above 69.96 Gy & 0.1 Gy & 0.00 & 0.00\\
Mean overdose SMG (R) above 69.96 Gy & 0.1 Gy & 0.00 & 0.00\\
\bottomrule
\end{tabular}
\end{table}
\begin{table}
\caption{Formulation and dose statistics for TROTS case 3}
\hskip-1.0cm\begin{tabular}{llll}
\toprule
Description & Bound / weight & IPOPT & Nymph \\
\midrule
Robust mean underdose CTV High below 64.68 Gy & 0.1 Gy & 0.10 & 0.10\\
Mean underdose CTV Intermediate 10 mm below 52.92 Gy & 0.1 Gy & 0.02 & 0.05\\
Mean underdose CTV Low Shrunk 10 mm below 52.92 Gy & 0.1 Gy & 0.10 & 0.10\\
Mean overdose Parotid (L) above 69.96 Gy & 0.1 Gy & 0.00 & 0.00\\
Mean overdose Parotid (R) above 69.96 Gy & 0.1 Gy & 0.00 & 0.00\\
Mean overdose SMG (L) above 69.96 Gy & 0.1 Gy & 0.00 & 0.00\\
Mean overdose SMG (R) above 69.96 Gy & 0.1 Gy & 0.01 & 0.00\\
\bottomrule
\end{tabular}
\end{table}
\begin{figure}
\centering
\includegraphics{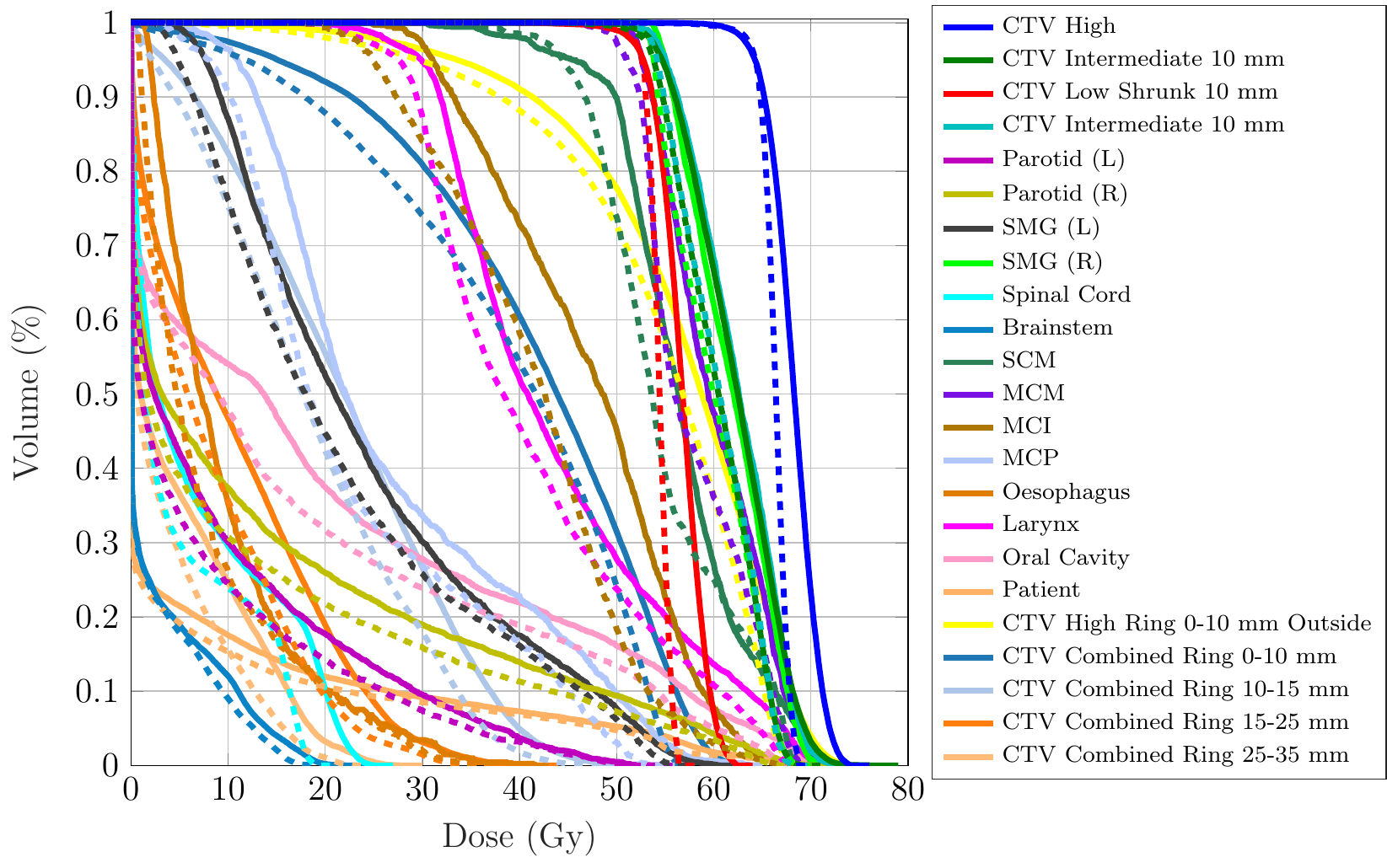}
\caption{DVH for TROTS case 3. Solid curves: L-BFGS. Dashed curves: Nymph.}
\end{figure}
\begin{figure}
\centering
\includegraphics{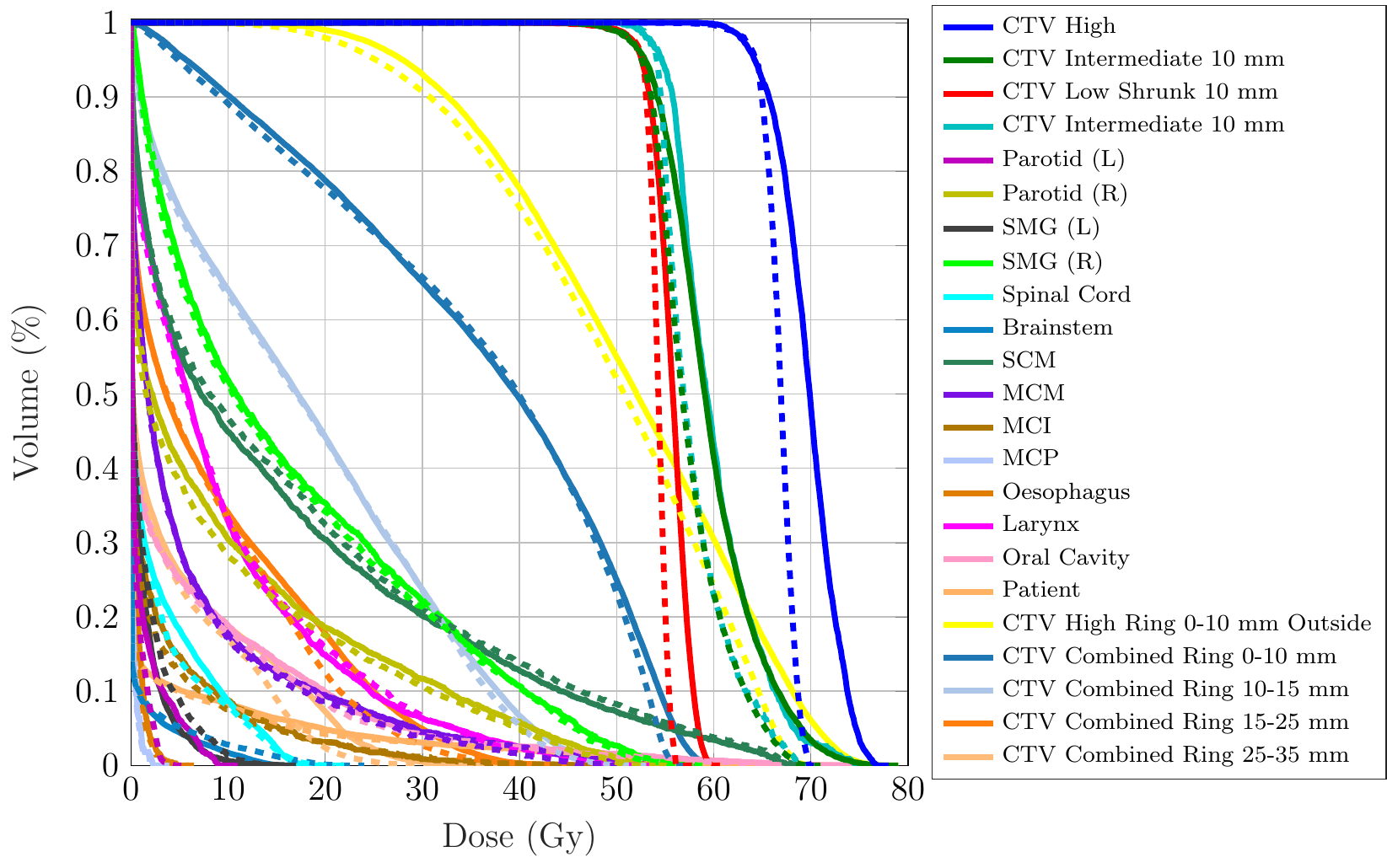}
\caption{DVH for TROTS case 4. Solid curves: L-BFGS. Dashed curves: Nymph.}
\end{figure}
\begin{table}
\caption{Formulation and dose statistics for TROTS case 4}
\hskip-1.0cm\begin{tabular}{llll}
\toprule
Description & Bound / weight & IPOPT & Nymph \\
\midrule
Robust mean underdose CTV High below 64.68 Gy & 0.1 Gy & 0.10 & 0.10\\
Mean underdose CTV Intermediate 10 mm below 52.92 Gy & 0.1 Gy & 0.10 & 0.10\\
Mean underdose CTV Low Shrunk 10 mm below 52.92 Gy & 0.1 Gy & 0.10 & 0.10\\
Mean overdose Parotid (L) above 69.96 Gy & 0.1 Gy & 0.00 & 0.00\\
Mean overdose Parotid (R) above 69.96 Gy & 0.1 Gy & 0.00 & 0.00\\
Mean overdose SMG (L) above 69.96 Gy & 0.1 Gy & 0.00 & 0.00\\
Mean overdose SMG (R) above 69.96 Gy & 0.1 Gy & 0.00 & 0.00\\
\bottomrule
\end{tabular}
\end{table}
\begin{table}
\caption{Formulation and dose statistics for TROTS case 5}
\hskip-1.0cm\begin{tabular}{llll}
\toprule
Description & Bound / weight & IPOPT & Nymph \\
\midrule
Robust mean underdose CTV High below 64.68 Gy & 0.1 Gy & 0.10 & 0.10\\
Mean underdose CTV Intermediate 10 mm below 52.92 Gy & 0.1 Gy & 0.10 & 0.09\\
Mean underdose CTV Low Shrunk 10 mm below 52.92 Gy & 0.1 Gy & 0.10 & 0.10\\
Mean overdose Parotid (L) above 69.96 Gy & 0.1 Gy & 0.00 & 0.00\\
Mean overdose Parotid (R) above 69.96 Gy & 0.1 Gy & 0.00 & 0.00\\
Mean overdose SMG (L) above 69.96 Gy & 0.1 Gy & 0.00 & 0.00\\
Mean overdose SMG (R) above 69.96 Gy & 0.1 Gy & 0.00 & 0.00\\
\bottomrule
\end{tabular}
\end{table}
\begin{figure}
\centering
\includegraphics{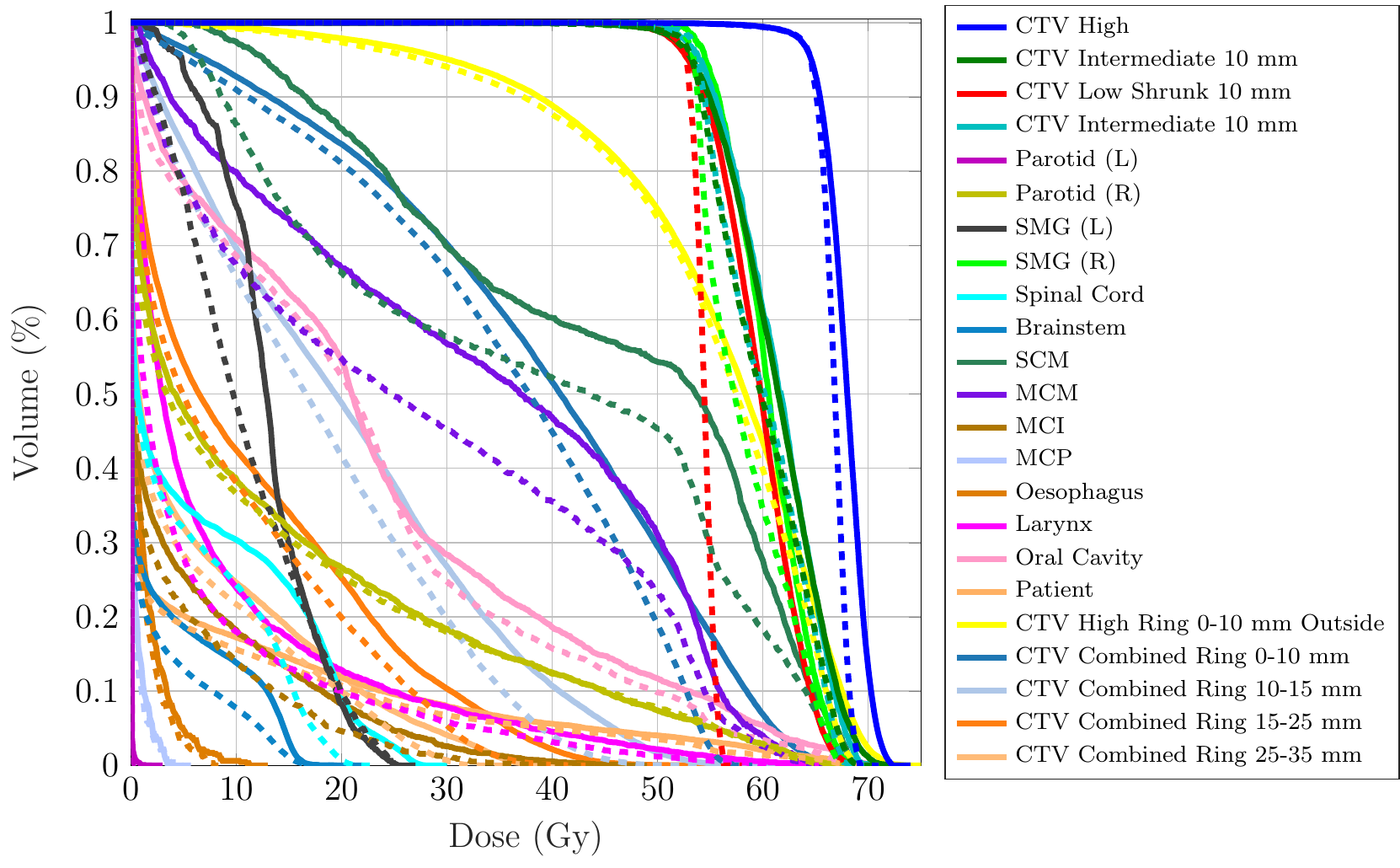}
\caption{DVH for TROTS case 5. Solid curves: L-BFGS. Dashed curves: Nymph.}
\end{figure}
\begin{figure}
\centering
\includegraphics{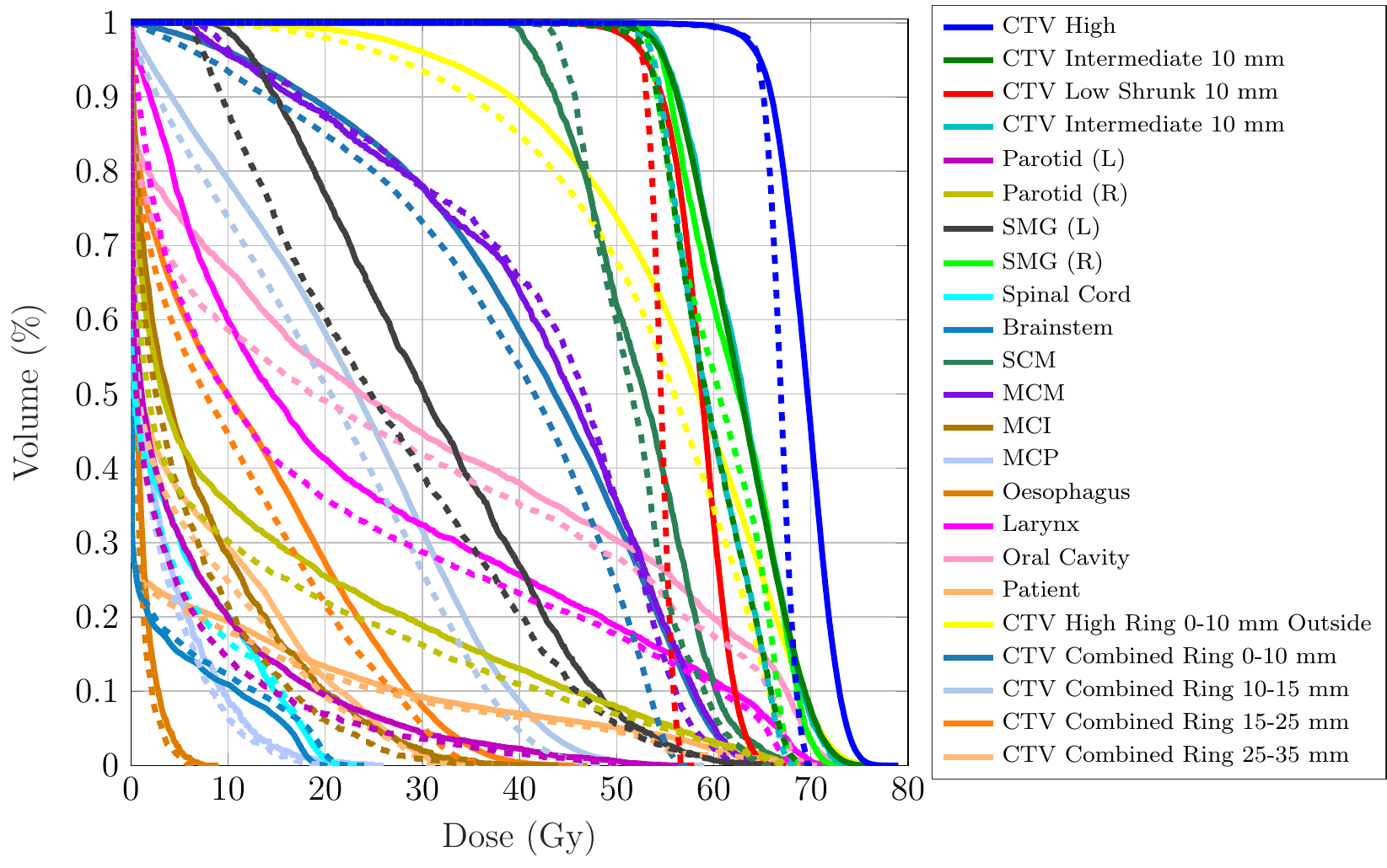}
\caption{DVH for TROTS case 6. Solid curves: L-BFGS. Dashed curves: Nymph.}
\end{figure}
\begin{table}
\caption{Formulation and dose statistics for TROTS case 6}
\hskip-1.0cm\begin{tabular}{llll}
\toprule
Description & Bound / weight & IPOPT & Nymph \\
\midrule
Robust mean underdose CTV High below 64.68 Gy & 0.1 Gy & 0.10 & 0.10\\
Mean underdose CTV Intermediate 10 mm below 52.92 Gy & 0.1 Gy & 0.03 & 0.08\\
Mean underdose CTV Low Shrunk 10 mm below 52.92 Gy & 0.1 Gy & 0.10 & 0.10\\
Mean overdose Parotid (L) above 69.96 Gy & 0.1 Gy & 0.00 & 0.00\\
Mean overdose Parotid (R) above 69.96 Gy & 0.1 Gy & 0.00 & 0.00\\
Mean overdose SMG (L) above 69.96 Gy & 0.1 Gy & 0.00 & 0.00\\
Mean overdose SMG (R) above 69.96 Gy & 0.1 Gy & 0.03 & 0.00\\
\bottomrule
\end{tabular}
\end{table}
\begin{table}
\caption{Formulation and dose statistics for TROTS case 7}
\hskip-1.0cm\begin{tabular}{llll}
\toprule
Description & Bound / weight & IPOPT & Nymph \\
\midrule
Robust mean underdose CTV High below 64.68 Gy & 0.1 Gy & 0.09 & 0.10\\
Mean underdose CTV Intermediate 10 mm below 52.92 Gy & 0.1 Gy & 0.05 & 0.08\\
Mean underdose CTV Low Shrunk 10 mm below 52.92 Gy & 0.1 Gy & 0.10 & 0.10\\
Mean overdose Parotid (L) above 69.96 Gy & 0.1 Gy & 0.00 & 0.00\\
Mean overdose Parotid (R) above 69.96 Gy & 0.1 Gy & 0.00 & 0.00\\
Mean overdose SMG (L) above 69.96 Gy & 0.1 Gy & 0.00 & 0.00\\
Mean overdose SMG (R) above 69.96 Gy & 0.1 Gy & 0.00 & 0.00\\
\bottomrule
\end{tabular}
\end{table}
\begin{figure}
\centering
\includegraphics{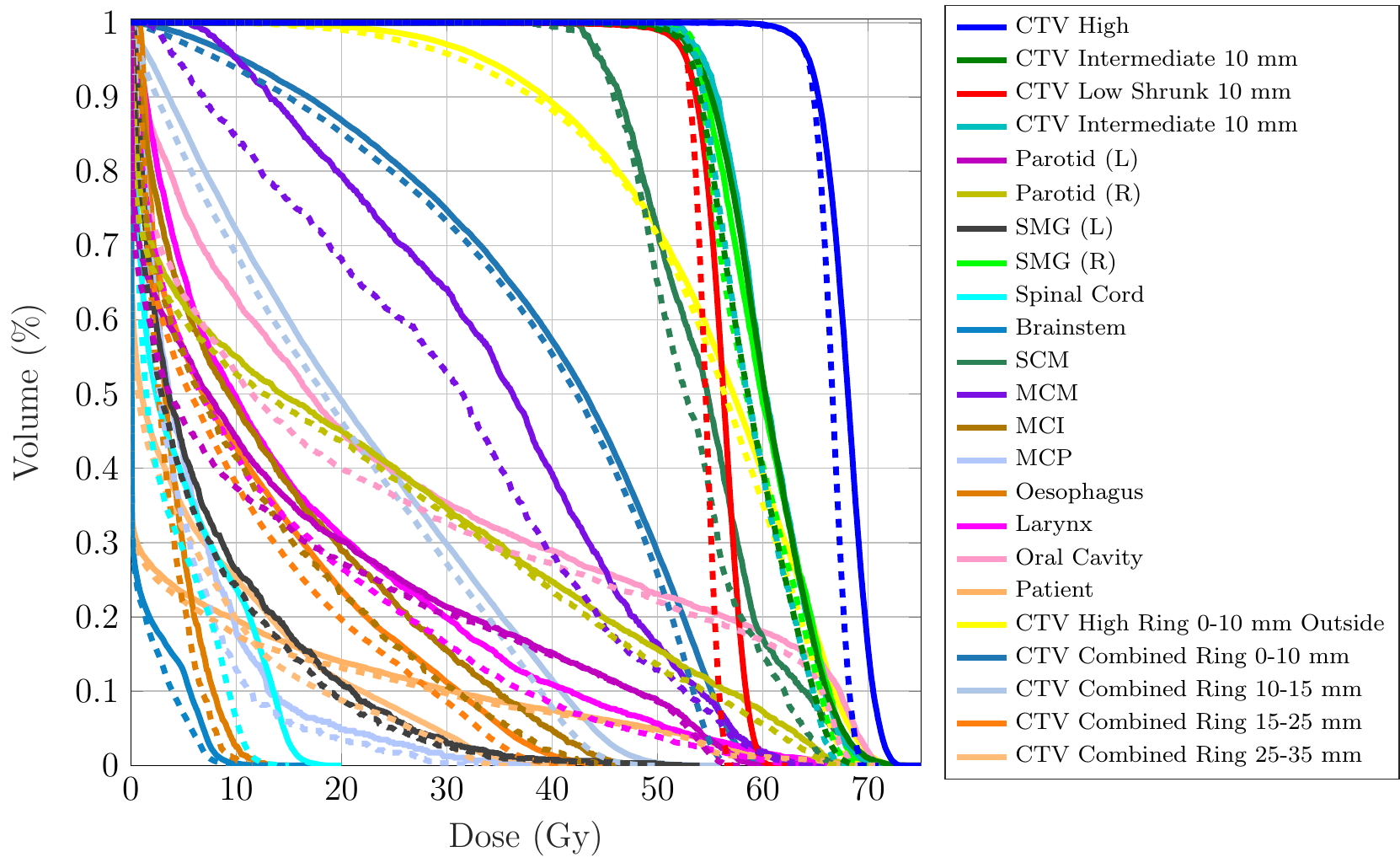}
\caption{DVH for TROTS case 7. Solid curves: L-BFGS. Dashed curves: Nymph.}
\end{figure}
\begin{figure}
\centering
\includegraphics{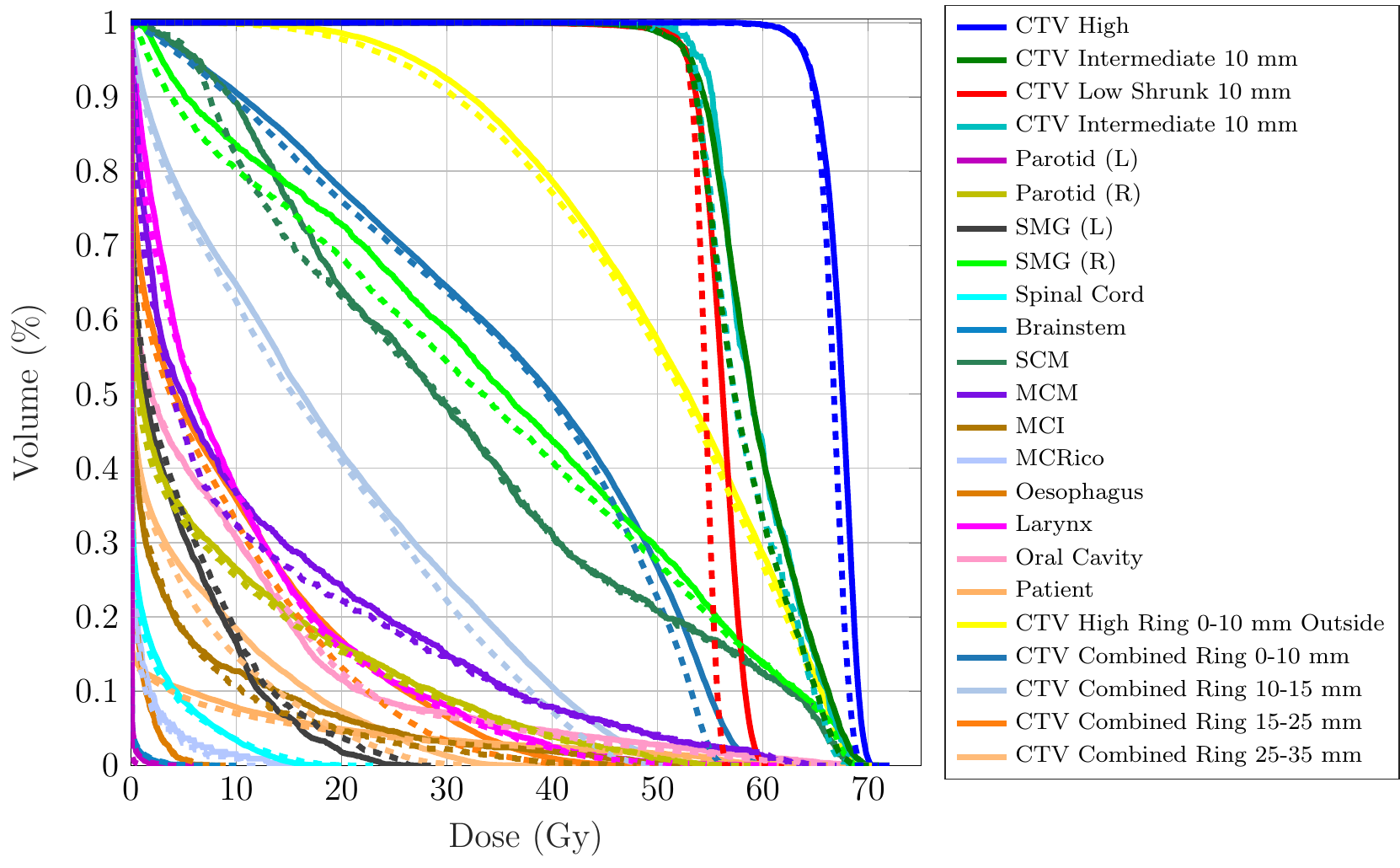}
\caption{DVH for TROTS case 8. Solid curves: L-BFGS. Dashed curves: Nymph.}
\end{figure}
\begin{table}
\caption{Formulation and dose statistics for TROTS case 8}
\hskip-1.0cm\begin{tabular}{llll}
\toprule
Description & Bound / weight & IPOPT & Nymph \\
\midrule
Robust mean underdose CTV High below 64.68 Gy & 0.1 Gy & 0.10 & 0.10\\
Mean underdose CTV Intermediate 10 mm below 52.92 Gy & 0.1 Gy & 0.10 & 0.10\\
Mean underdose CTV Low Shrunk 10 mm below 52.92 Gy & 0.1 Gy & 0.10 & 0.10\\
Mean overdose Parotid (L) above 69.96 Gy & 0.1 Gy & 0.00 & 0.00\\
Mean overdose Parotid (R) above 69.96 Gy & 0.1 Gy & 0.00 & 0.00\\
Mean overdose SMG (L) above 69.96 Gy & 0.1 Gy & 0.00 & 0.00\\
Mean overdose SMG (R) above 69.96 Gy & 0.1 Gy & 0.00 & 0.00\\
\bottomrule
\end{tabular}
\end{table}
\begin{table}
\caption{Formulation and dose statistics for TROTS case 9}
\hskip-1.0cm\begin{tabular}{llll}
\toprule
Description & Bound / weight & IPOPT & Nymph \\
\midrule
Robust mean underdose CTV High below 64.68 Gy & 0.1 Gy & 0.10 & 0.10\\
Mean underdose CTV Intermediate 10 mm below 52.92 Gy & 0.1 Gy & 0.06 & 0.08\\
Mean underdose CTV Low Shrunk 10 mm below 52.92 Gy & 0.1 Gy & 0.12 & 0.10\\
Mean overdose Parotid (L) above 69.96 Gy & 0.1 Gy & 0.04 & 0.00\\
Mean overdose Parotid (R) above 69.96 Gy & 0.1 Gy & 0.00 & 0.00\\
Mean overdose SMG (L) above 69.96 Gy & 0.1 Gy & 0.08 & 0.00\\
Mean overdose SMG (R) above 69.96 Gy & 0.1 Gy & 0.00 & 0.00\\
\bottomrule
\end{tabular}
\end{table}
\begin{figure}
\centering
\includegraphics{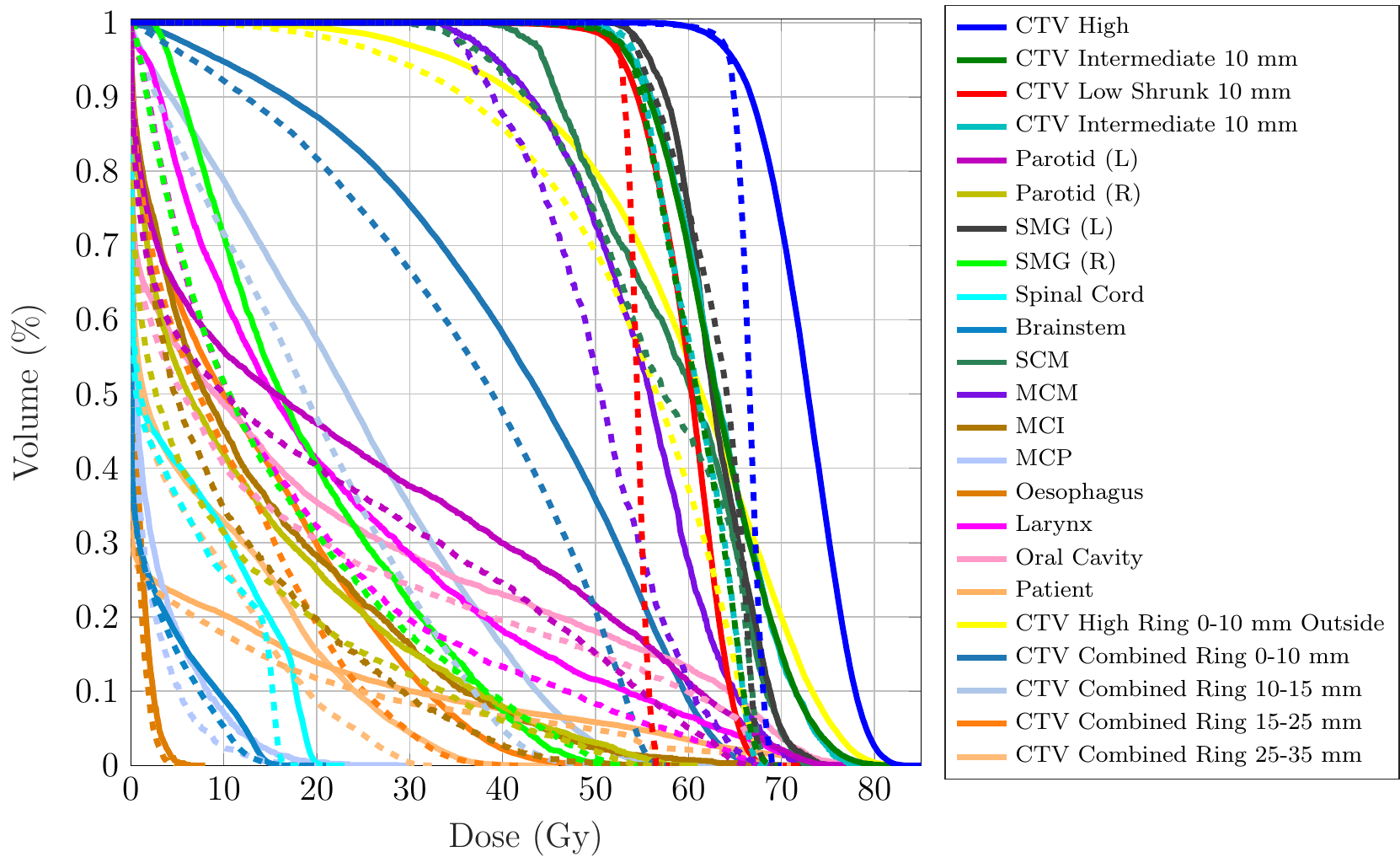}
\caption{DVH for TROTS case 9. Solid curves: L-BFGS. Dashed curves: Nymph.}
\end{figure}
\begin{figure}
\centering
\includegraphics{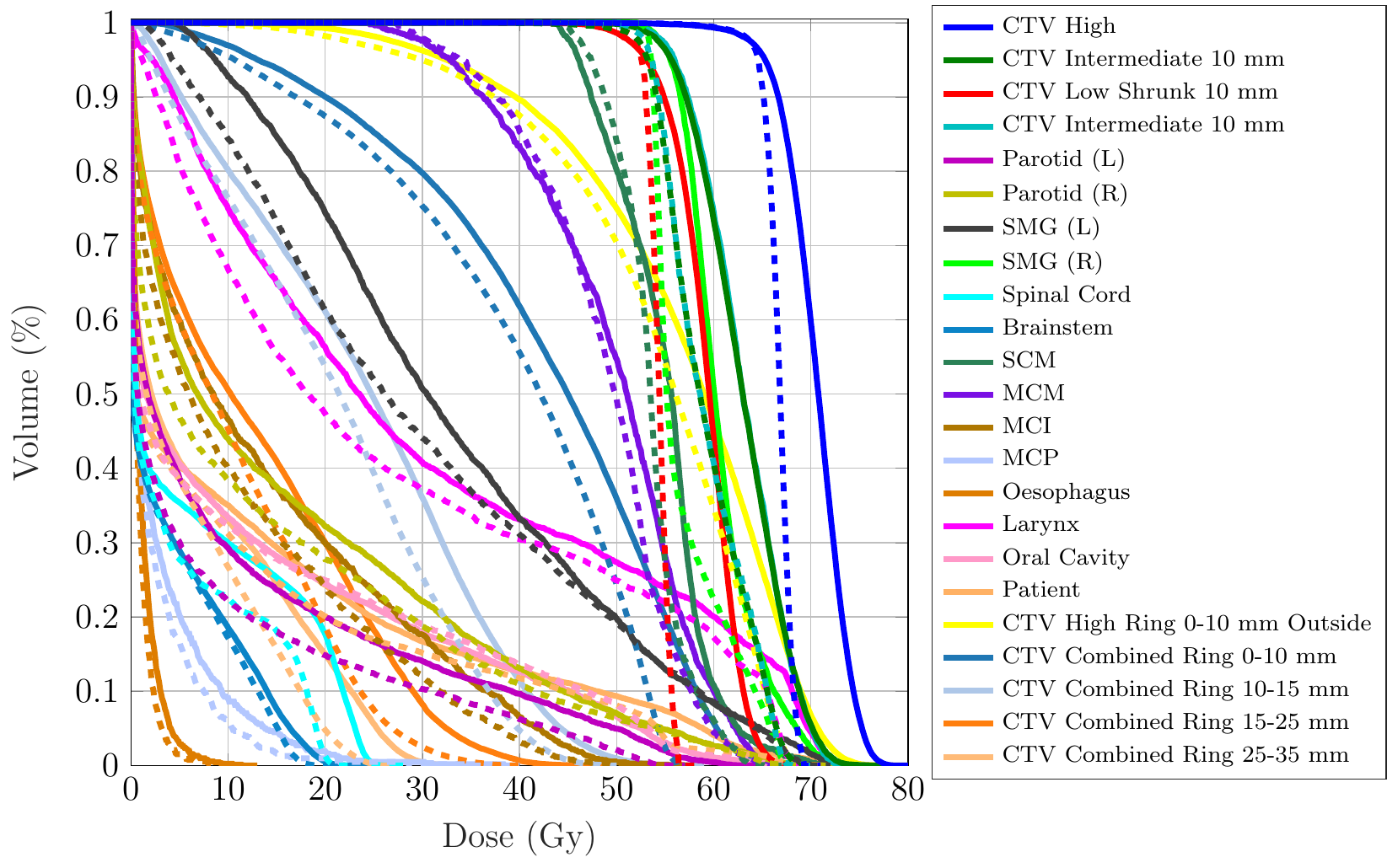}
\caption{DVH for TROTS case 10. Solid curves: L-BFGS. Dashed curves: Nymph.}
\end{figure}
\begin{table}
\caption{Formulation and dose statistics for TROTS case 10}
\hskip-1.0cm\begin{tabular}{llll}
\toprule
Description & Bound / weight & IPOPT & Nymph \\
\midrule
Robust mean underdose CTV High below 64.68 Gy & 0.1 Gy & 0.10 & 0.10\\
Mean underdose CTV Intermediate 10 mm below 52.92 Gy & 0.1 Gy & 0.02 & 0.08\\
Mean underdose CTV Low Shrunk 10 mm below 52.92 Gy & 0.1 Gy & 0.11 & 0.10\\
Mean overdose Parotid (L) above 69.96 Gy & 0.1 Gy & 0.00 & 0.00\\
Mean overdose Parotid (R) above 69.96 Gy & 0.1 Gy & 0.00 & 0.00\\
Mean overdose SMG (L) above 69.96 Gy & 0.1 Gy & 0.03 & 0.00\\
Mean overdose SMG (R) above 69.96 Gy & 0.1 Gy & 0.03 & 0.00\\
\bottomrule
\end{tabular}
\end{table}
\begin{table}
\caption{Formulation and dose statistics for TROTS case 11}
\hskip-1.0cm\begin{tabular}{llll}
\toprule
Description & Bound / weight & IPOPT & Nymph \\
\midrule
Robust mean underdose CTV High below 64.68 Gy & 0.1 Gy & 0.12 & 0.10\\
Mean underdose CTV Intermediate 10 mm below 52.92 Gy & 0.1 Gy & 0.11 & 0.10\\
Mean underdose CTV Low Shrunk 10 mm below 52.92 Gy & 0.1 Gy & 0.16 & 0.10\\
Mean overdose Parotid (L) above 69.96 Gy & 0.1 Gy & 0.00 & 0.00\\
Mean overdose Parotid (R) above 69.96 Gy & 0.1 Gy & 0.00 & 0.00\\
Mean overdose SMG (L) above 69.96 Gy & 0.1 Gy & 0.00 & 0.00\\
Mean overdose SMG (R) above 69.96 Gy & 0.1 Gy & 0.00 & 0.00\\
\bottomrule
\end{tabular}
\end{table}
\begin{figure}
\centering
\includegraphics{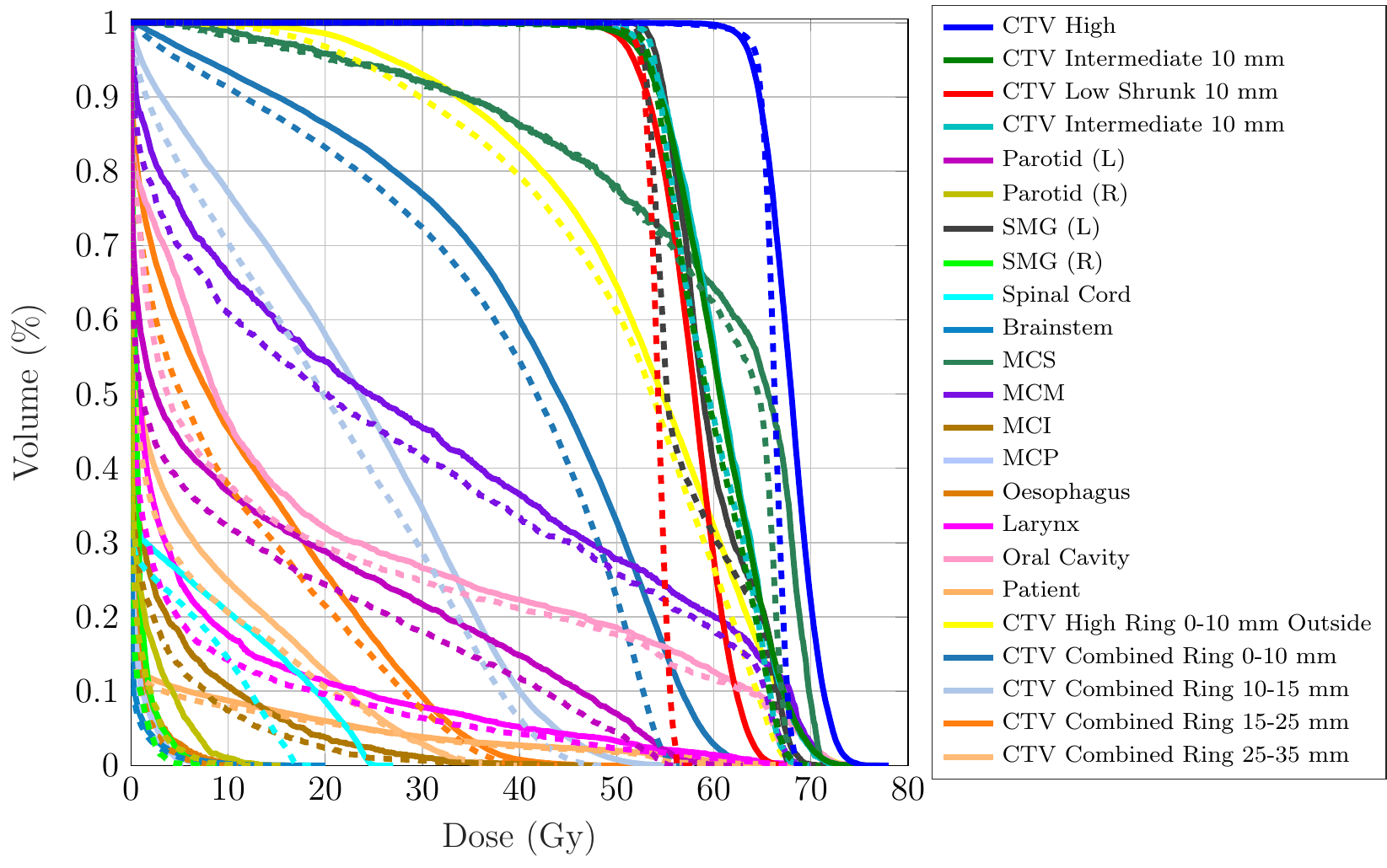}
\caption{DVH for TROTS case 11. Solid curves: L-BFGS. Dashed curves: Nymph.}
\end{figure}
\begin{figure}
\centering
\includegraphics{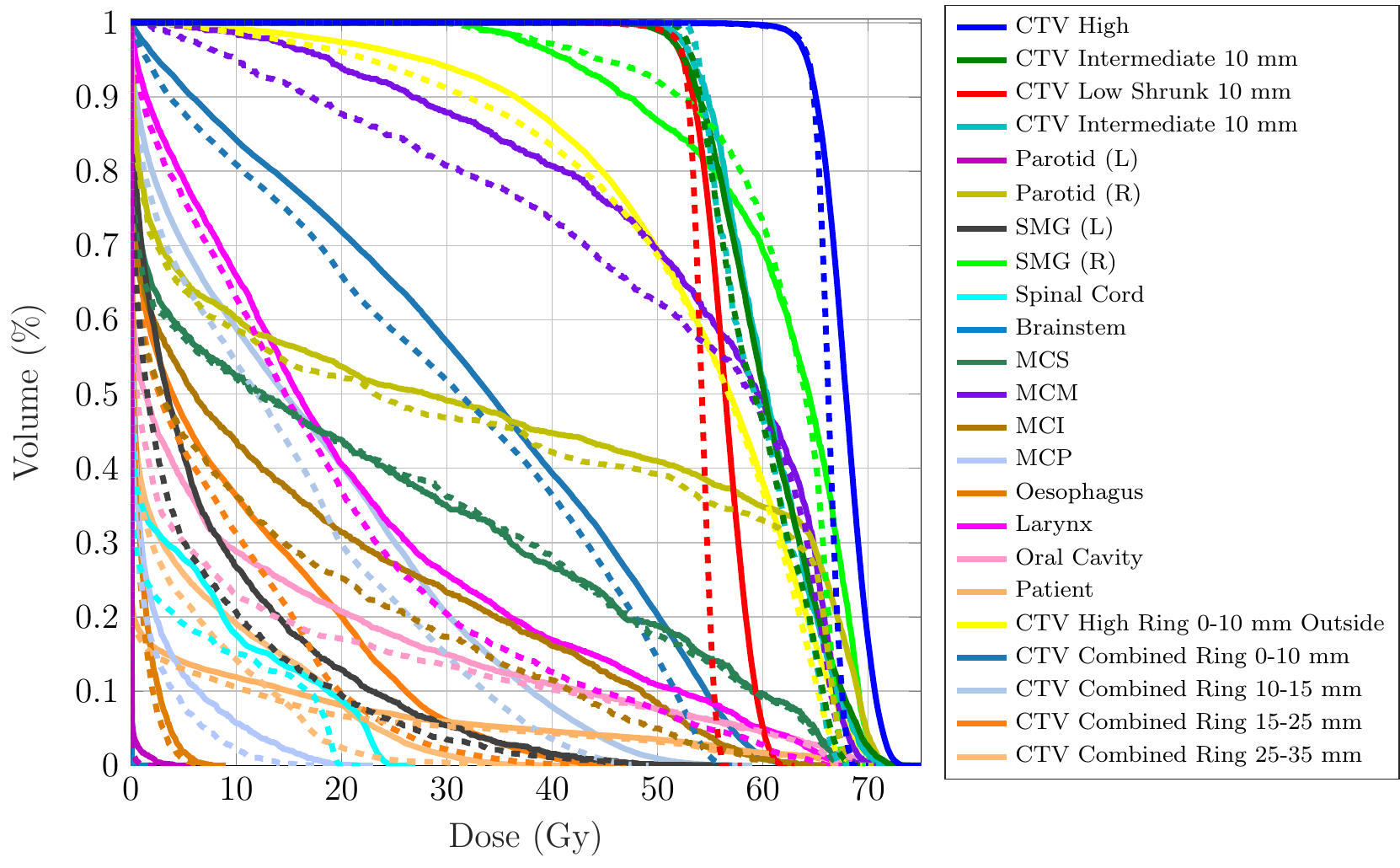}
\caption{DVH for TROTS case 12. Solid curves: L-BFGS. Dashed curves: Nymph.}
\end{figure}
\begin{table}
\caption{Formulation and dose statistics for TROTS case 12}
\hskip-1.0cm\begin{tabular}{llll}
\toprule
Description & Bound / weight & IPOPT & Nymph \\
\midrule
Robust mean underdose CTV High below 64.68 Gy & 0.1 Gy & 0.10 & 0.10\\
Mean underdose CTV Intermediate 10 mm below 52.92 Gy & 0.1 Gy & 0.09 & 0.02\\
Mean underdose CTV Low Shrunk 10 mm below 52.92 Gy & 0.1 Gy & 0.10 & 0.10\\
Mean overdose Parotid (L) above 69.96 Gy & 0.1 Gy & 0.00 & 0.00\\
Mean overdose Parotid (R) above 69.96 Gy & 0.1 Gy & 0.03 & 0.00\\
Mean overdose SMG (L) above 69.96 Gy & 0.1 Gy & 0.00 & 0.00\\
Mean overdose SMG (R) above 69.96 Gy & 0.1 Gy & 0.03 & 0.00\\
\bottomrule
\end{tabular}
\end{table}
\begin{table}
\caption{Formulation and dose statistics for TROTS case 13}
\hskip-1.0cm\begin{tabular}{llll}
\toprule
Description & Bound / weight & IPOPT & Nymph \\
\midrule
Robust mean underdose CTV High below 64.68 Gy & 0.1 Gy & 0.11 & 0.10\\
Mean underdose CTV Intermediate 10 mm below 52.92 Gy & 0.1 Gy & 0.04 & 0.06\\
Mean underdose CTV Low Shrunk 10 mm below 52.92 Gy & 0.1 Gy & 0.14 & 0.10\\
Mean overdose Parotid (L) above 69.96 Gy & 0.1 Gy & 0.00 & 0.00\\
Mean overdose Parotid (R) above 69.96 Gy & 0.1 Gy & 0.00 & 0.00\\
Mean overdose SMG (L) above 69.96 Gy & 0.1 Gy & 0.04 & 0.00\\
Mean overdose SMG (R) above 69.96 Gy & 0.1 Gy & 0.00 & 0.00\\
\bottomrule
\end{tabular}
\end{table}
\begin{figure}
\centering
\includegraphics{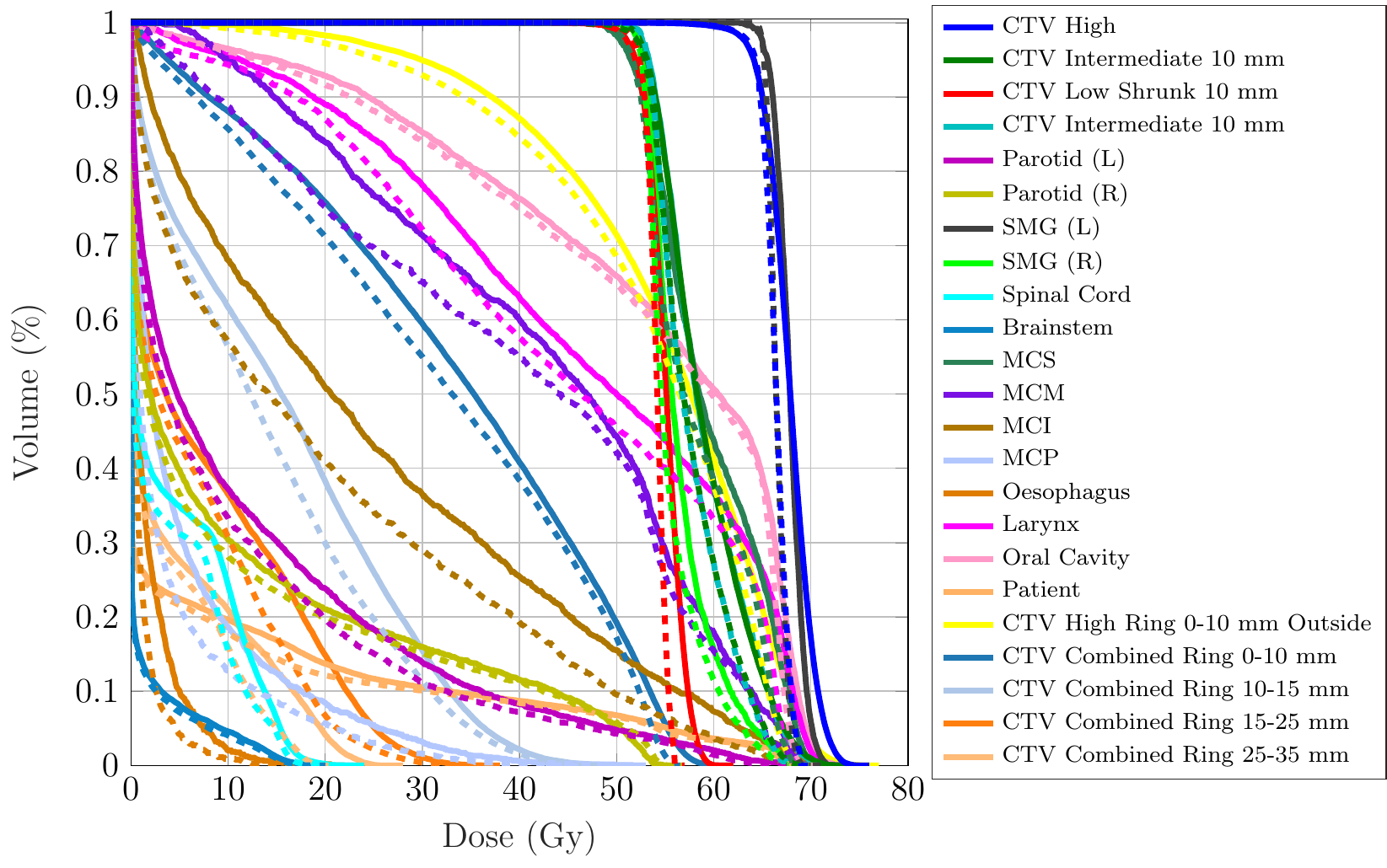}
\caption{DVH for TROTS case 13. Solid curves: L-BFGS. Dashed curves: Nymph.}
\end{figure}
\begin{figure}
\centering
\includegraphics{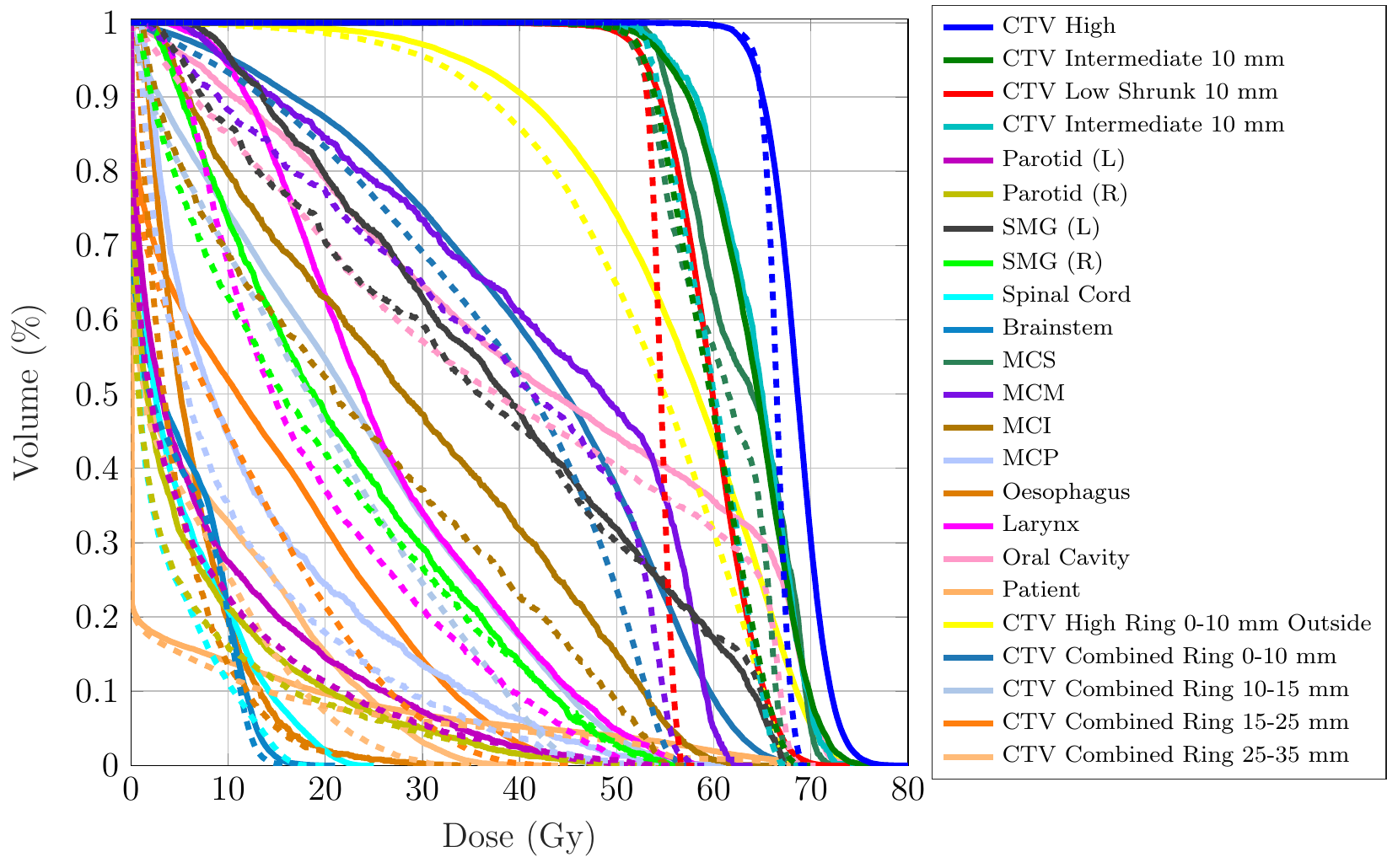}
\caption{DVH for TROTS case 14. Solid curves: L-BFGS. Dashed curves: Nymph.}
\end{figure}
\begin{table}
\caption{Formulation and dose statistics for TROTS case 14}
\hskip-1.0cm\begin{tabular}{llll}
\toprule
Description & Bound / weight & IPOPT & Nymph \\
\midrule
Robust mean underdose CTV High below 64.68 Gy & 0.1 Gy & 0.11 & 0.10\\
Mean underdose CTV Intermediate 10 mm below 52.92 Gy & 0.1 Gy & 0.05 & 0.09\\
Mean underdose CTV Low Shrunk 10 mm below 52.92 Gy & 0.1 Gy & 0.12 & 0.10\\
Mean overdose Parotid (L) above 69.96 Gy & 0.1 Gy & 0.00 & 0.00\\
Mean overdose Parotid (R) above 69.96 Gy & 0.1 Gy & 0.00 & 0.00\\
Mean overdose SMG (L) above 69.96 Gy & 0.1 Gy & 0.00 & 0.00\\
Mean overdose SMG (R) above 69.96 Gy & 0.1 Gy & 0.00 & 0.00\\
\bottomrule
\end{tabular}
\end{table}
\begin{table}
\caption{Formulation and dose statistics for TROTS case 15}
\hskip-1.0cm\begin{tabular}{llll}
\toprule
Description & Bound / weight & IPOPT & Nymph \\
\midrule
Robust mean underdose CTV High below 64.68 Gy & 0.1 Gy & 0.10 & 0.10\\
Mean underdose CTV Intermediate 10 mm below 52.92 Gy & 0.1 Gy & 0.08 & 0.10\\
Mean underdose CTV Low Shrunk 10 mm below 52.92 Gy & 0.1 Gy & 0.09 & 0.10\\
Mean overdose Parotid (L) above 69.96 Gy & 0.1 Gy & 0.00 & 0.00\\
Mean overdose Parotid (R) above 69.96 Gy & 0.1 Gy & 0.00 & 0.00\\
Mean overdose SMG (L) above 69.96 Gy & 0.1 Gy & 0.00 & 0.00\\
Mean overdose SMG (R) above 69.96 Gy & 0.1 Gy & 0.00 & 0.00\\
\bottomrule
\end{tabular}
\end{table}
\begin{figure}
\centering
\includegraphics{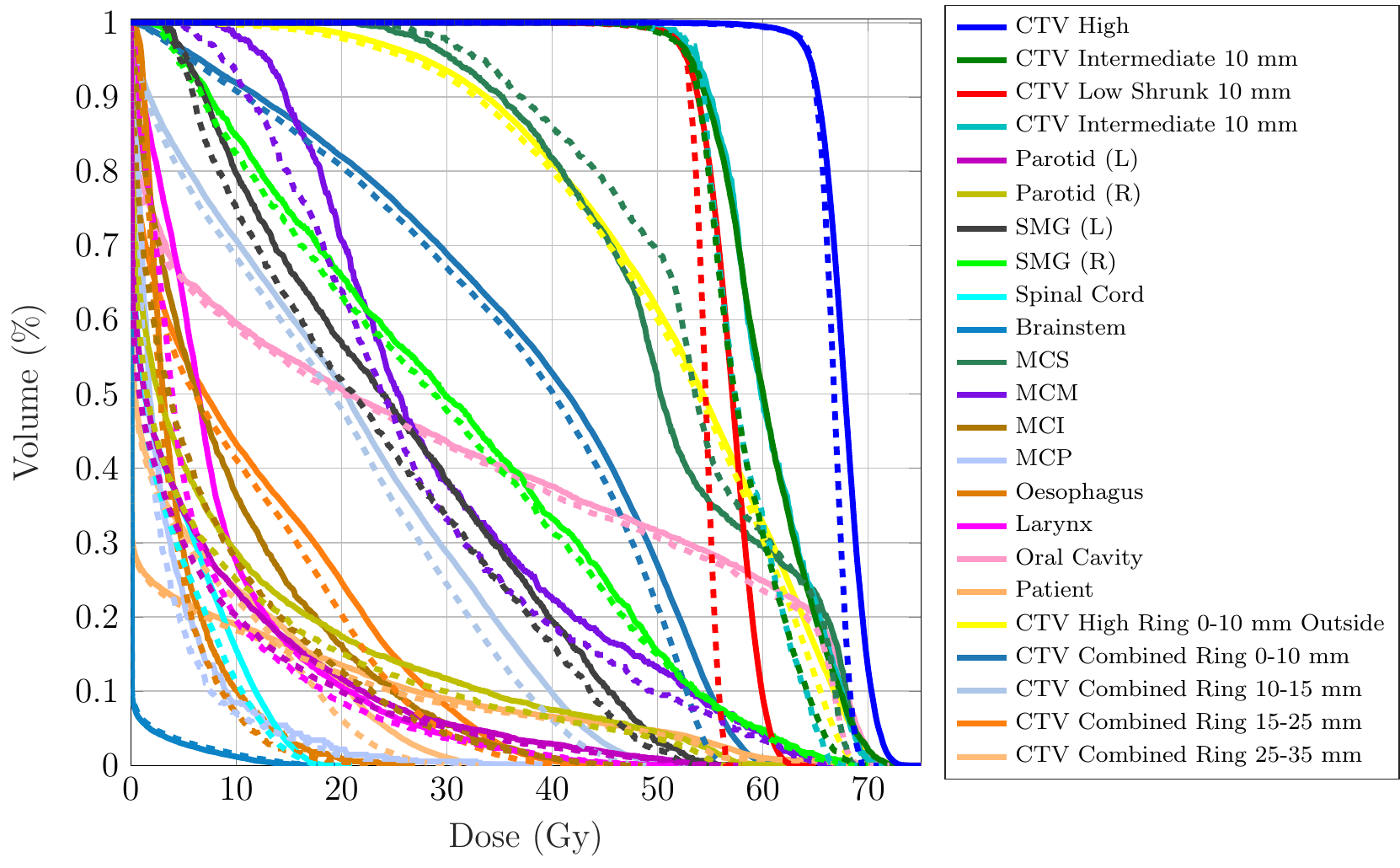}
\caption{DVH for TROTS case 15. Solid curves: L-BFGS. Dashed curves: Nymph.}
\end{figure}
\begin{figure}
\centering
\includegraphics{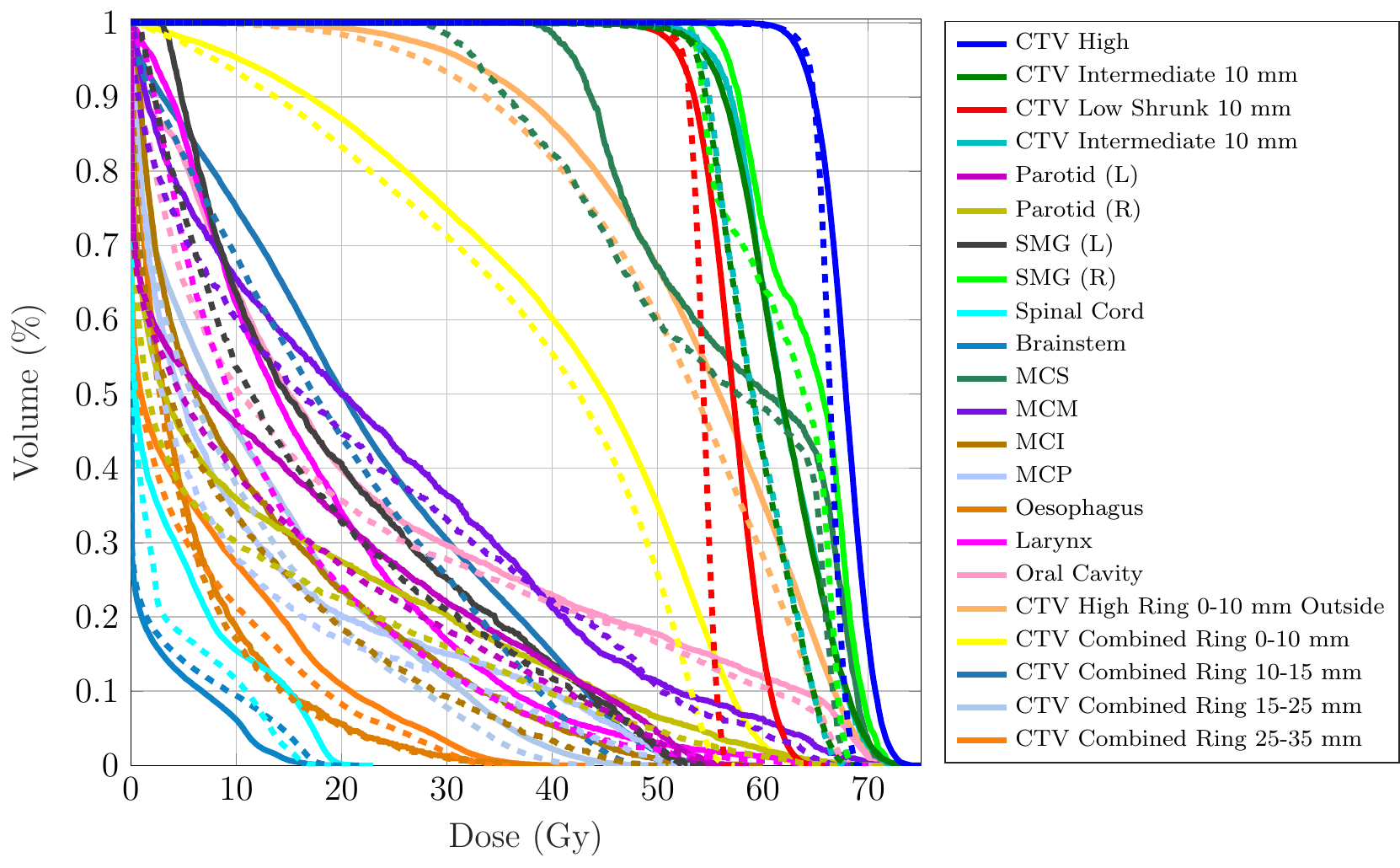}
\caption{DVH for TROTS case 16. Solid curves: L-BFGS. Dashed curves: Nymph.}
\end{figure}
\begin{table}
\caption{Formulation and dose statistics for TROTS case 16}
\hskip-1.0cm\begin{tabular}{llll}
\toprule
Description & Bound / weight & IPOPT & Nymph \\
\midrule
Robust mean underdose CTV High below 64.68 Gy & 0.1 Gy & 0.12 & 0.10\\
Mean underdose CTV Intermediate 10 mm below 52.92 Gy & 0.1 Gy & 0.05 & 0.07\\
Mean underdose CTV Low Shrunk 10 mm below 52.92 Gy & 0.1 Gy & 0.12 & 0.10\\
Mean overdose Parotid (L) above 69.96 Gy & 0.1 Gy & 0.00 & 0.00\\
Mean overdose Parotid (R) above 69.96 Gy & 0.1 Gy & 0.00 & 0.00\\
Mean overdose SMG (L) above 69.96 Gy & 0.1 Gy & 0.00 & 0.00\\
Mean overdose SMG (R) above 69.96 Gy & 0.1 Gy & 0.06 & 0.00\\
\bottomrule
\end{tabular}
\end{table}
\begin{table}
\caption{Formulation and dose statistics for TROTS case 17}
\hskip-1.0cm\begin{tabular}{llll}
\toprule
Description & Bound / weight & IPOPT & Nymph \\
\midrule
Robust mean underdose CTV High below 64.68 Gy & 0.1 Gy & 0.10 & 0.10\\
Mean underdose CTV Intermediate 10 mm below 52.92 Gy & 0.1 Gy & 0.10 & 0.10\\
Mean underdose CTV Low Shrunk 10 mm below 52.92 Gy & 0.1 Gy & 0.10 & 0.09\\
Mean overdose Parotid (L) above 69.96 Gy & 0.1 Gy & 0.00 & 0.00\\
Mean overdose Parotid (R) above 69.96 Gy & 0.1 Gy & 0.00 & 0.00\\
Mean overdose SMG (L) above 69.96 Gy & 0.1 Gy & 0.00 & 0.00\\
Mean overdose SMG (R) above 69.96 Gy & 0.1 Gy & 0.00 & 0.00\\
\bottomrule
\end{tabular}
\end{table}
\begin{figure}
\centering
\includegraphics{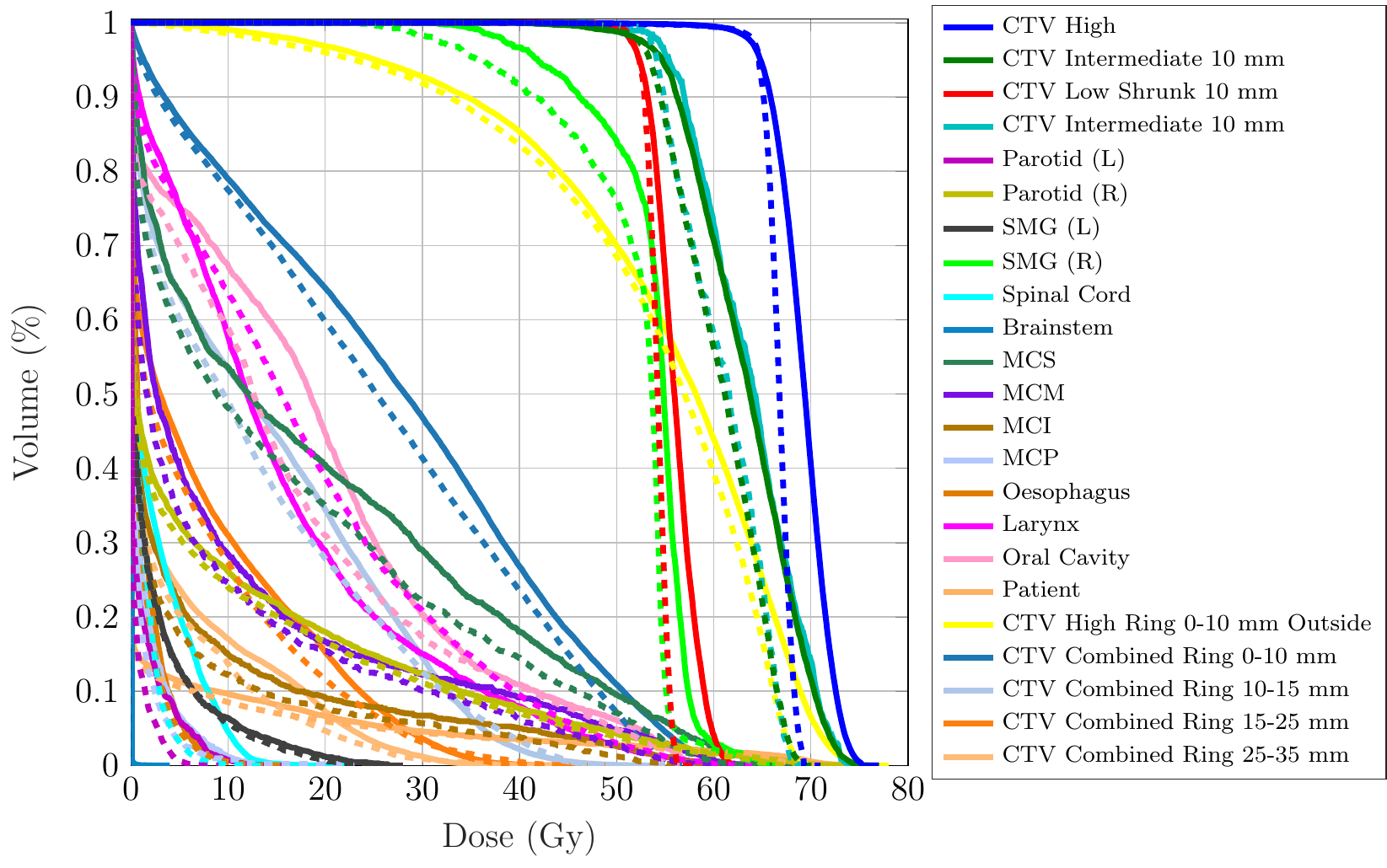}
\caption{DVH for TROTS case 17. Solid curves: L-BFGS. Dashed curves: Nymph.}
\end{figure}
\begin{figure}
\centering
\includegraphics{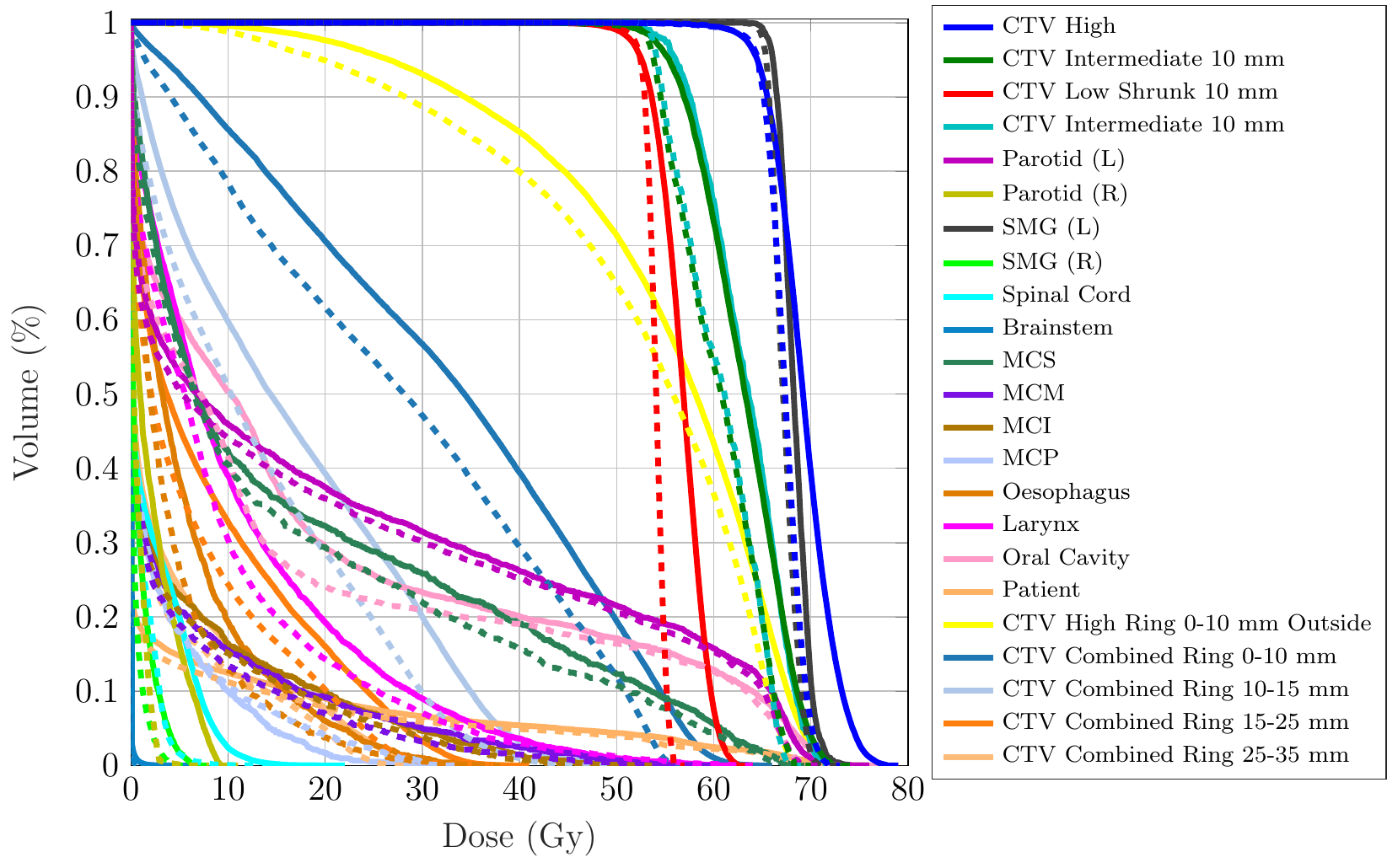}
\caption{DVH for TROTS case 18. Solid curves: L-BFGS. Dashed curves: Nymph.}
\end{figure}
\begin{table}
\caption{Formulation and dose statistics for TROTS case 18}
\hskip-1.0cm\begin{tabular}{llll}
\toprule
Description & Bound / weight & IPOPT & Nymph \\
\midrule
Robust mean underdose CTV High below 64.68 Gy & 0.1 Gy & 0.10 & 0.10\\
Mean underdose CTV Intermediate 10 mm below 52.92 Gy & 0.1 Gy & 0.03 & 0.02\\
Mean underdose CTV Low Shrunk 10 mm below 52.92 Gy & 0.1 Gy & 0.10 & 0.09\\
Mean overdose Parotid (L) above 69.96 Gy & 0.1 Gy & 0.01 & 0.00\\
Mean overdose Parotid (R) above 69.96 Gy & 0.1 Gy & 0.00 & 0.00\\
Mean overdose SMG (L) above 69.96 Gy & 0.1 Gy & 0.10 & 0.01\\
Mean overdose SMG (R) above 69.96 Gy & 0.1 Gy & 0.00 & 0.00\\
\bottomrule
\end{tabular}
\end{table}
\begin{table}
\caption{Formulation and dose statistics for TROTS case 19}
\hskip-1.0cm\begin{tabular}{llll}
\toprule
Description & Bound / weight & IPOPT & Nymph \\
\midrule
Robust mean underdose CTV High below 64.68 Gy & 0.1 Gy & 0.31 & 0.09\\
Mean underdose CTV Intermediate 10 mm below 52.92 Gy & 0.1 Gy & 0.03 & 0.08\\
Mean underdose CTV Low Shrunk 10 mm below 52.92 Gy & 0.1 Gy & 0.21 & 0.10\\
Mean overdose Parotid (L) above 69.96 Gy & 0.1 Gy & 0.00 & 0.00\\
Mean overdose Parotid (R) above 69.96 Gy & 0.1 Gy & 0.00 & 0.00\\
Mean overdose SMG (L) above 69.96 Gy & 0.1 Gy & 0.03 & 0.00\\
Mean overdose SMG (R) above 69.96 Gy & 0.1 Gy & 0.00 & 0.00\\
\bottomrule
\end{tabular}
\end{table}
\begin{figure}
\centering
\includegraphics{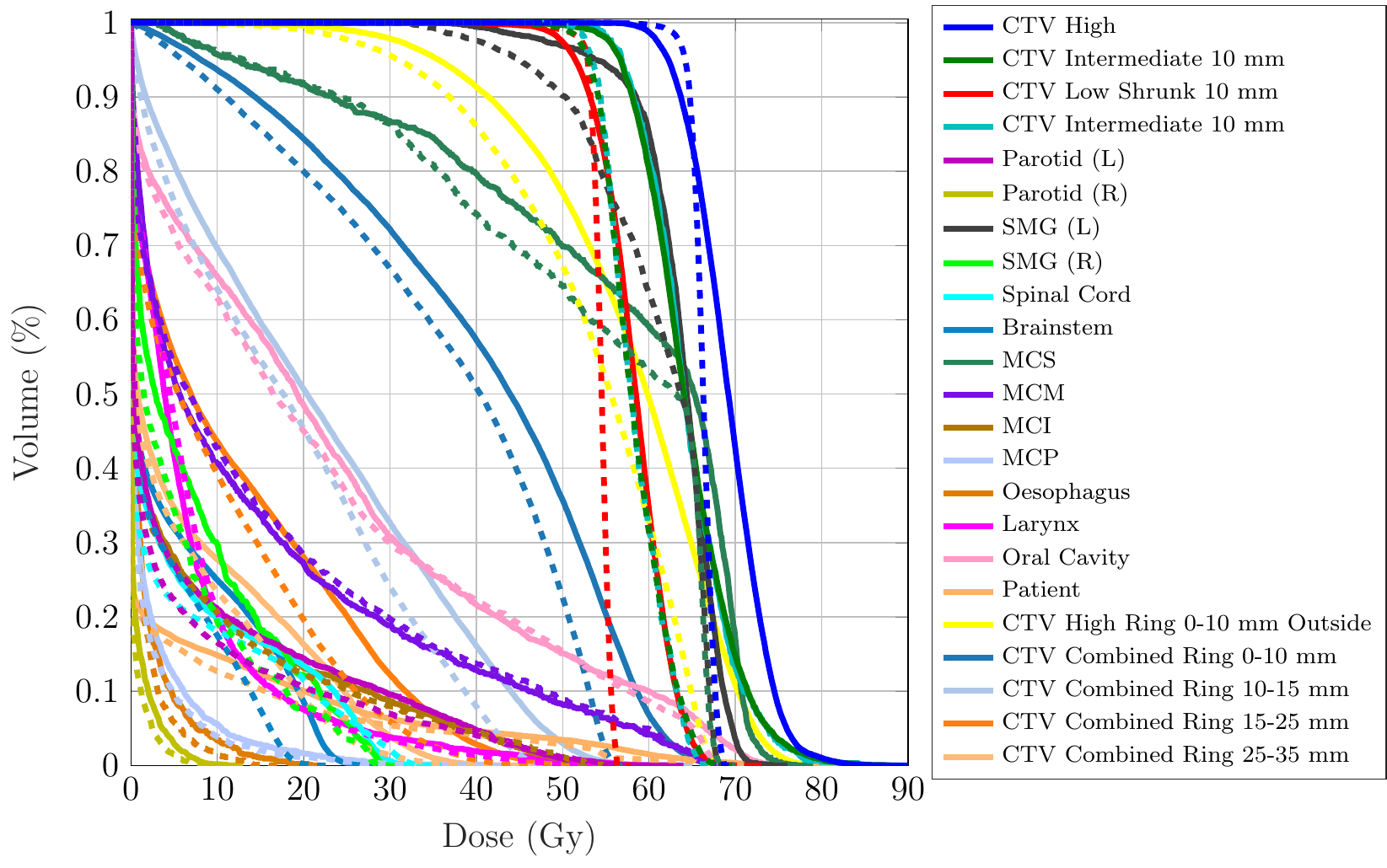}
\caption{DVH for TROTS case 19. Solid curves: L-BFGS. Dashed curves: Nymph.}
\end{figure}
\begin{figure}
\centering
\includegraphics{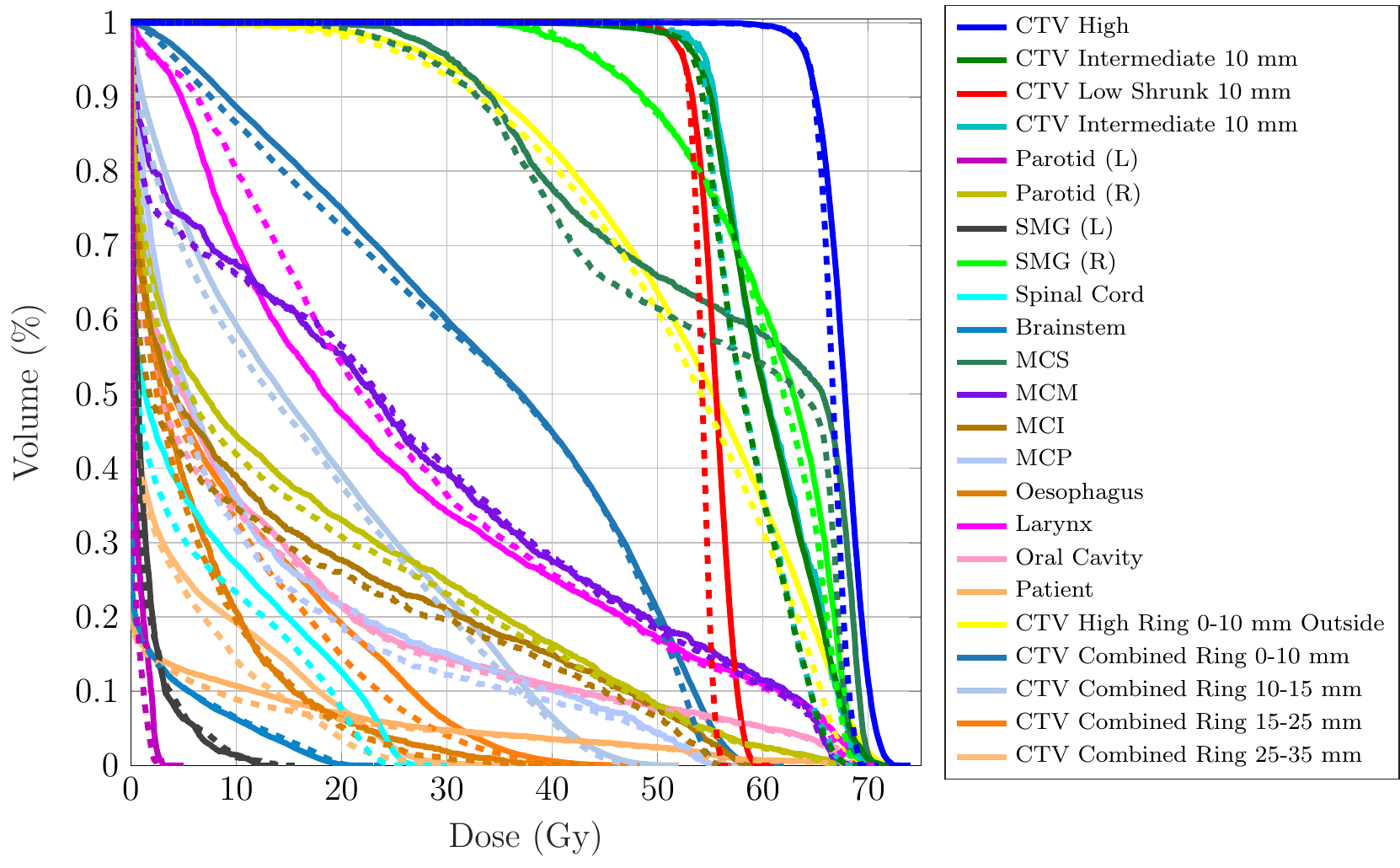}
\caption{DVH for TROTS case 20. Solid curves: L-BFGS. Dashed curves: Nymph.}
\end{figure}
\begin{table}
\caption{Formulation and dose statistics for TROTS case 20}
\hskip-1.0cm\begin{tabular}{llll}
\toprule
Description & Bound / weight & IPOPT & Nymph \\
\midrule
Robust mean underdose CTV High below 64.68 Gy & 0.1 Gy & 0.10 & 0.10\\
Mean underdose CTV Intermediate 10 mm below 52.92 Gy & 0.1 Gy & 0.10 & 0.10\\
Mean underdose CTV Low Shrunk 10 mm below 52.92 Gy & 0.1 Gy & 0.10 & 0.10\\
Mean overdose Parotid (L) above 69.96 Gy & 0.1 Gy & 0.00 & 0.00\\
Mean overdose Parotid (R) above 69.96 Gy & 0.1 Gy & 0.00 & 0.00\\
Mean overdose SMG (L) above 69.96 Gy & 0.1 Gy & 0.00 & 0.00\\
Mean overdose SMG (R) above 69.96 Gy & 0.1 Gy & 0.00 & 0.00\\
\bottomrule
\end{tabular}
\end{table}


\begin{thebibliography}{47}
\providecommand{\natexlab}[1]{#1}
\providecommand{\url}[1]{\texttt{#1}}
\expandafter\ifx\csname urlstyle\endcsname\relax
  \providecommand{\doi}[1]{doi: #1}\else
  \providecommand{\doi}{doi: \begingroup \urlstyle{rm}\Url}\fi

\bibitem[Alterovitz et~al.(2006)Alterovitz, Lessard, Pouliot, Hsu, O'Brien, and
  Goldberg]{alterovitz2006optimization}
R.~Alterovitz, E.~Lessard, J.~Pouliot, I.-C.~J. Hsu, J.~F. O'Brien, and
  K.~Goldberg.
\newblock Optimization of {HDR} brachytherapy dose distributions using linear
  programming with penalty costs.
\newblock \emph{Medical Physics}, 33\penalty0 (11):\penalty0 4012--4019, 2006.

\bibitem[Andersen and Andersen(2000)]{andersen2000mosek}
E.~D. Andersen and K.~D. Andersen.
\newblock The {MOSEK} interior point optimizer for linear programming: an
  implementation of the homogeneous algorithm.
\newblock In \emph{High performance optimization}, 197--232. Springer, 2000.

\bibitem[Andersen et~al.(2011)Andersen, Dahl, Liu, and
  Vandenberghe]{andersen2011interior}
M.~Andersen, J.~Dahl, Z.~Liu, and L.~Vandenberghe.
\newblock \emph{Interior-point methods for large-scale cone programming}.
\newblock MIT Press, Cambridge, MA, USA, 2011.

\bibitem[Bertsimas and Sim(2004)]{bertsimas2004price}
D.~Bertsimas and M.~Sim.
\newblock The price of robustness.
\newblock \emph{Operations Research}, 52\penalty0 (1):\penalty0 35--53, 2004.

\bibitem[Bixby(2012)]{bixby2012brief}
R.~E. Bixby.
\newblock A brief history of linear and mixed-integer programming computation.
\newblock \emph{Documenta Mathematica}, 107--121, 2012.

\bibitem[Bortfeld et~al.(1990)Bortfeld, B{\"u}rkelbach, Boesecke, and
  Schlegel]{bortfeld1990methods}
T.~Bortfeld, J.~B{\"u}rkelbach, R.~Boesecke, and W.~Schlegel.
\newblock Methods of image reconstruction from projections applied to
  conformation radiotherapy.
\newblock \emph{Physics in Medicine \& Biology}, 35\penalty0 (10):\penalty0
  1423, 1990.

\bibitem[Boyd and Vandenberghe(2004)]{Boydcvx}
S.~Boyd and L.~Vandenberghe.
\newblock \emph{Convex Optimization}.
\newblock Cambridge University Press, New York, NY, USA, 2004.

\bibitem[Breedveld and Heijmen(2017)]{breedveld2017data}
S.~Breedveld and B.~Heijmen.
\newblock Data for {TROTS}--the radiotherapy optimisation test set.
\newblock \emph{Data in brief}, 12:\penalty0 143--149, 2017.

\bibitem[Breedveld et~al.(2017)Breedveld, van~den Berg, and
  Heijmen]{Breedveld2016}
S.~Breedveld, B.~van~den Berg, and B.~Heijmen.
\newblock An interior-point implementation developed and tuned for radiation
  therapy treatment planning.
\newblock \emph{Computational Optimization and Applications}, 68\penalty0
  (2):\penalty0 209--242, 2017.

\bibitem[Breedveld et~al.(2019)Breedveld, Craft, Van~Haveren, and
  Heijmen]{breedveld2019multi}
S.~Breedveld, D.~Craft, R.~Van~Haveren, and B.~Heijmen.
\newblock Multi-criteria optimization and decision-making in radiotherapy.
\newblock \emph{European Journal of Operational Research}, 277\penalty0
  (1):\penalty0 1--19, 2019.

\bibitem[Carlsson et~al.(2006)Carlsson, Forsgren, Rehbinder, and
  Eriksson]{carlsson2006using}
F.~Carlsson, A.~Forsgren, H.~Rehbinder, and K.~Eriksson.
\newblock Using eigenstructure of the {H}essian to reduce the dimension of the
  intensity modulated radiation therapy optimization problem.
\newblock \emph{Annals of Operations Research}, 148\penalty0 (1):\penalty0
  81--94, 2006.

\bibitem[Castro(2016)]{castro2016interior}
J.~Castro.
\newblock Interior-point solver for convex separable block-angular problems.
\newblock \emph{Optimization Methods and Software}, 31\penalty0 (1):\penalty0
  88--109, 2016.

\bibitem[Charnes et~al.(1955)Charnes, Cooper, and Ferguson]{charnes1955optimal}
A.~Charnes, W.~W. Cooper, and R.~O. Ferguson.
\newblock Optimal estimation of executive compensation by linear programming.
\newblock \emph{Management Science}, 1\penalty0 (2):\penalty0 138--151, 1955.

\bibitem[Chen et~al.(2010)Chen, Craft, Madden, Zhang, Kooy, and
  Herman]{chen2010fast}
W.~Chen, D.~Craft, T.~M. Madden, K.~Zhang, H.~M. Kooy, and G.~T. Herman.
\newblock A fast optimization algorithm for multicriteria intensity modulated
  proton therapy planning.
\newblock \emph{Medical Physics}, 37\penalty0 (9):\penalty0 4938--4945, 2010.

\bibitem[Colombo and Gondzio(2008)]{colombo2008further}
M.~Colombo and J.~Gondzio.
\newblock Further development of multiple centrality correctors for interior
  point methods.
\newblock \emph{Computational Optimization and Applications}, 41\penalty0
  (3):\penalty0 277--305, 2008.

\bibitem[Craft and Bortfeld(2008)]{craft2008many}
D.~Craft and T.~Bortfeld.
\newblock How many plans are needed in an {IMRT} multi-objective plan database?
\newblock \emph{Physics in Medicine \& Biology}, 53\penalty0 (11):\penalty0
  2785, 2008.

\bibitem[Engberg et~al.(2017)Engberg, Forsgren, Eriksson, and
  H{\aa}rdemark]{engberg2017explicit}
L.~Engberg, A.~Forsgren, K.~Eriksson, and B.~H{\aa}rdemark.
\newblock Explicit optimization of plan quality measures in intensity-modulated
  radiation therapy treatment planning.
\newblock \emph{Medical Physics}, 44\penalty0 (6):\penalty0 2045--2053, 2017.

\bibitem[Fourer et~al.(2002)Fourer, Gay, and Kernighan]{fourer1997ampl}
R.~Fourer, D.~Gay, and B.~Kernighan.
\newblock \emph{AMPL: a modeling language for mathematical programming}.
\newblock Duxbury Press, 2002.

\bibitem[Fragniere et~al.(2000)Fragniere, Gondzio, Sarkissian, and
  Vial]{fragniere2000structure}
E.~Fragniere, J.~Gondzio, R.~Sarkissian, and J.-P. Vial.
\newblock A structure-exploiting tool in algebraic modeling languages.
\newblock \emph{Management Science}, 46\penalty0 (8):\penalty0 1145--1158,
  2000.

\bibitem[Fredriksson and Bokrantz(2014)]{fredriksson2014critical}
A.~Fredriksson and R.~Bokrantz.
\newblock A critical evaluation of worst case optimization methods for robust
  intensity-modulated proton therapy planning.
\newblock \emph{Medical Physics}, 41\penalty0 (8):\penalty0 081701, 2014.

\bibitem[Gay(1985)]{gay1985electronic}
D.~M. Gay.
\newblock Electronic mail distribution of linear programming test problems.
\newblock \emph{Mathematical Programming Society COAL Newsletter}, 13:\penalty0
  10--12, 1985.

\bibitem[Gill et~al.(2005)Gill, Murray, and Saunders]{gill2005snopt}
P.~E. Gill, W.~Murray, and M.~A. Saunders.
\newblock {SNOPT}: An {SQP} algorithm for large-scale constrained optimization.
\newblock \emph{SIAM Review}, 47\penalty0 (1):\penalty0 99--131, 2005.

\bibitem[Gondzio(1996)]{gondzio1996multiple}
J.~Gondzio.
\newblock Multiple centrality corrections in a primal-dual method for linear
  programming.
\newblock \emph{Computational Optimization and Applications}, 6\penalty0
  (2):\penalty0 137--156, 1996.

\bibitem[Gondzio and Grothey(2009)]{gondzio2009exploiting}
J.~Gondzio and A.~Grothey.
\newblock Exploiting structure in parallel implementation of interior point
  methods for optimization.
\newblock \emph{Computational Management Science}, 6\penalty0 (2):\penalty0
  135--160, 2009.

\bibitem[Gondzio and Sarkissian(2003)]{gondzio2003parallel}
J.~Gondzio and R.~Sarkissian.
\newblock Parallel interior-point solver for structured linear programs.
\newblock \emph{Mathematical Programming}, 96\penalty0 (3):\penalty0 561--584,
  2003.

\bibitem[Gondzio et~al.(1997)Gondzio, Sarkissian, and Vial]{gondzio1997using}
J.~Gondzio, R.~Sarkissian, and J.-P. Vial.
\newblock Using an interior point method for the master problem in a
  decomposition approach.
\newblock \emph{European Journal of Operational Research}, 101\penalty0
  (3):\penalty0 577--587, 1997.

\bibitem[Gorissen(2020)]{Gorissen2020}
B.~L. Gorissen.
\newblock Nymph, the fastest exact inverse planning algorithm for radiation
  therapy, 2020.
\newblock URL \url{https://3142.nl/nymph/}.

\bibitem[Gorissen et~al.(2013)Gorissen, Den~Hertog, and
  Hoffmann]{gorissen2013mixed}
B.~L. Gorissen, D.~Den~Hertog, and A.~L. Hoffmann.
\newblock Mixed integer programming improves comprehensibility and plan quality
  in inverse optimization of prostate {HDR} brachytherapy.
\newblock \emph{Physics in Medicine \& Biology}, 58\penalty0 (4):\penalty0
  1041--1058, 2013.

\bibitem[Guttman(1946)]{guttman1946enlargement}
L.~Guttman.
\newblock Enlargement methods for computing the inverse matrix.
\newblock \emph{The annals of mathematical statistics}, 17\penalty0
  (3):\penalty0 336--343, 1946.

\bibitem[Karabis et~al.(2009)Karabis, Belotti, and
  Baltas]{karabis2009optimization}
A.~Karabis, P.~Belotti, and D.~Baltas.
\newblock Optimization of catheter position and dwell time in prostate {HDR}
  brachytherapy using {HIPO} and linear programming.
\newblock In \emph{World Congress on Medical Physics and Biomedical
  Engineering, September 7-12, 2009, Munich, Germany}, 612--615, 2009.

\bibitem[Karp(1972)]{karp1972reducibility}
R.~M. Karp.
\newblock Reducibility among combinatorial problems.
\newblock In \emph{Complexity of computer computations}, 85--103. Springer,
  1972.

\bibitem[Lessard and Pouliot(2001)]{lessard2001inverse}
E.~Lessard and J.~Pouliot.
\newblock Inverse planning anatomy-based dose optimization for
  {HDR}-brachytherapy of the prostate using fast simulated annealing algorithm
  and dedicated objective function.
\newblock \emph{Medical Physics}, 28\penalty0 (5):\penalty0 773--779, 2001.

\bibitem[Lomax(1999)]{lomax1999intensity}
A.~Lomax.
\newblock Intensity modulation methods for proton radiotherapy.
\newblock \emph{Physics in Medicine \& Biology}, 44\penalty0 (1):\penalty0
  185--205, 1999.

\bibitem[Mehrotra(1992)]{mehrotra1992implementation}
S.~Mehrotra.
\newblock On the implementation of a primal-dual interior point method.
\newblock \emph{SIAM Journal on Optimization}, 2\penalty0 (4):\penalty0
  575--601, 1992.

\bibitem[Mukherjee et~al.(2020)Mukherjee, Hong, Deasy, and
  Zarepisheh]{mukherjee2019integrating}
S.~Mukherjee, L.~Hong, J.~O. Deasy, and M.~Zarepisheh.
\newblock Integrating soft and hard dose-volume constraints into hierarchical
  constrained {IMRT} optimization.
\newblock \emph{Medical Physics}, 47\penalty0 (2):\penalty0 414--421, 2020.

\bibitem[Rockafellar and Uryasev(2000)]{rockafellar2000optimization}
R.~T. Rockafellar and S.~Uryasev.
\newblock Optimization of conditional value-at-risk.
\newblock \emph{Journal of risk}, 2\penalty0 (3):\penalty0 21--42, 2000.

\bibitem[Romeijn et~al.(2003)Romeijn, Ahuja, Dempsey, Kumar, and
  Li]{romeijn2003novel}
H.~E. Romeijn, R.~K. Ahuja, J.~F. Dempsey, A.~Kumar, and J.~G. Li.
\newblock A novel linear programming approach to fluence map optimization for
  intensity modulated radiation therapy treatment planning.
\newblock \emph{Physics in Medicine \& Biology}, 48\penalty0 (21):\penalty0
  3521--3542, 2003.

\bibitem[Romeijn et~al.(2006)Romeijn, Ahuja, Dempsey, and
  Kumar]{romeijn2006new}
H.~E. Romeijn, R.~K. Ahuja, J.~F. Dempsey, and A.~Kumar.
\newblock A new linear programming approach to radiation therapy treatment
  planning problems.
\newblock \emph{Operations Research}, 54\penalty0 (2):\penalty0 201--216, 2006.

\bibitem[Shepard et~al.(1999)Shepard, Ferris, Olivera, and
  Mackie]{shepard1999optimizing}
D.~M. Shepard, M.~C. Ferris, G.~H. Olivera, and T.~R. Mackie.
\newblock Optimizing the delivery of radiation therapy to cancer patients.
\newblock \emph{SIAM Review}, 41\penalty0 (4):\penalty0 721--744, 1999.

\bibitem[Ungun et~al.(2019)Ungun, Xing, and Boyd]{ungun2019real}
B.~Ungun, L.~Xing, and S.~Boyd.
\newblock Real-time radiation treatment planning with optimality guarantees via
  cluster and bound methods.
\newblock \emph{INFORMS Journal on Computing}, 31\penalty0 (3):\penalty0
  544--558, 2019.

\bibitem[Unkelbach et~al.(2018)Unkelbach, Alber, Bangert, Bokrantz, Chan,
  Deasy, Fredriksson, Gorissen, Van~Herk, Liu, Mahmoudzadeh, Nohadani, Siebers,
  Witte, and Xu]{unkelbach2018robust}
J.~Unkelbach, M.~Alber, M.~Bangert, R.~Bokrantz, T.~C.~Y. Chan, J.~O. Deasy,
  A.~Fredriksson, B.~L. Gorissen, M.~Van~Herk, W.~Liu, H.~Mahmoudzadeh,
  O.~Nohadani, J.~V. Siebers, M.~Witte, and H.~Xu.
\newblock Robust radiotherapy planning.
\newblock \emph{Physics in Medicine \& Biology}, 63\penalty0 (22):\penalty0
  22TR02, 2018.

\bibitem[Vanderbei(1999)]{vanderbei1999loqo}
R.~J. Vanderbei.
\newblock {LOQO}: An interior point code for quadratic programming.
\newblock \emph{Optimization methods and software}, 11\penalty0
  (1--4):\penalty0 451--484, 1999.

\bibitem[W{\"a}chter and Biegler(2006)]{wachter2006implementation}
A.~W{\"a}chter and L.~T. Biegler.
\newblock On the implementation of an interior-point filter line-search
  algorithm for large-scale nonlinear programming.
\newblock \emph{Mathematical Programming}, 106\penalty0 (1):\penalty0 25--57,
  2006.

\bibitem[Webb(2003)]{webb2003physical}
S.~Webb.
\newblock The physical basis of {IMRT} and inverse planning.
\newblock \emph{The British journal of radiology}, 76\penalty0 (910):\penalty0
  678--689, 2003.

\bibitem[Zhang(2017)]{zhang2017robust}
R.~Y. Zhang.
\newblock \emph{Robust stability analysis for large-scale power systems}.
\newblock PhD thesis, Massachusetts Institute of Technology, 2017.

\bibitem[Zhu et~al.(1997)Zhu, Byrd, Lu, and Nocedal]{zhu1997algorithm}
C.~Zhu, R.~H. Byrd, P.~Lu, and J.~Nocedal.
\newblock Algorithm 778: {L-BFGS-B}: {F}ortran subroutines for large-scale
  bound-constrained optimization.
\newblock \emph{ACM Transactions on mathematical software}, 23\penalty0
  (4):\penalty0 550--560, 1997.

\bibitem[Ziegenhein et~al.(2013)Ziegenhein, Kamerling, Bangert, Kunkel, and
  Oelfke]{ziegenhein2013performance}
P.~Ziegenhein, C.~P. Kamerling, M.~Bangert, J.~Kunkel, and U.~Oelfke.
\newblock Performance-optimized clinical {IMRT} planning on modern {CPU}s.
\newblock \emph{Physics in Medicine \& Biology}, 58\penalty0 (11):\penalty0
  3705--3715, 2013.

\end{thebibliography}
\end{document}